\pdfoutput=1
\documentclass[pdflatex,sn-mathphys-num]{sn-jnl}

\usepackage[T1]{fontenc}
\usepackage{lmodern}
\usepackage{microtype}
\usepackage{graphicx}
\usepackage{amsmath,amssymb,amsfonts}
\usepackage{amsthm}
\usepackage{enumitem}
\usepackage{xcolor}
\usepackage{etoolbox}
\makeatletter
\patchcmd{\@maketitle}{Contributing authors:}{E-mail:}{}{
  \PackageWarning{cmp-format}{Could not replace the author-email label}}
\makeatother

\hypersetup{colorlinks=true,linkcolor=blue,citecolor=blue,urlcolor=blue}

\numberwithin{equation}{section}

\theoremstyle{thmstyleone}
\newtheorem{theorem}{Theorem}[section]
\newtheorem{lemma}{Lemma}[section]
\newtheorem{proposition}{Proposition}[section]
\newtheorem{corollary}{Corollary}[section]

\theoremstyle{thmstylethree}

\newtheorem{assumption}{Assumption}[section]

\theoremstyle{thmstyletwo}

\newcommand{\bu}{\mathbf u}
\newcommand{\E}{\mathbb E}
\newcommand{\Var}{\operatorname{Var}}
\newcommand{\tr}{\operatorname{tr}}
\newcommand{\cyc}{\#\operatorname{cyc}}
\newcommand{\W}{\mathbf W}

\setlist[itemize]{leftmargin=*,itemsep=3pt,topsep=4pt}
\setlist[enumerate]{leftmargin=*,itemsep=3pt,topsep=4pt}
\begin{document}

\title[Large Powers of Ginibre Matrices]{Large Powers of Ginibre Matrices: Singular-Value Statistics via Eigenvalue Traces}

\author[1]{\fnm{Mariia} \sur{Seleznova}}\email{mariia.seleznova@tu-darmstadt.de}
\author[2]{\fnm{Nina} \sur{D\"ornemann}}\email{ndoernemann@math.au.dk}

\affil[1]{\orgdiv{Signal Theory and Statistical Learning Group}, \orgname{TU Darmstadt}, \orgaddress{\country{Germany}}}
\affil[2]{\orgdiv{Department of Mathematics}, \orgname{Aarhus University}, \orgaddress{\country{Denmark}}}

\abstract{\normalfont\unboldmath Repeated multiplication by a single nonnormal $n\times n$ matrix $\W$
governs many recurrent systems, including recurrent neural networks,
state-space models, and linear dynamical systems. Signal amplification and
stability in such systems depend on the singular values of the powers
$\W^k$, yet theory for singular values of large powers of random matrices
remains comparatively sparse. We study the average squared singular value
$Z_{k,n} = \frac1n\sum_{j=1}^n s_j(\W^k)^2$ 
for normalized complex Ginibre $\W$ when $k=k(n)$ grows with the dimension $n$.   
This joint limit lies between two classical asymptotic settings: fixed $n$
and growing $k$, common in dynamical systems, and fixed $k$ with growing $n$,
common in random matrix theory. Although
$Z_{k,n}$ has no evident samplewise relation to the eigenvalues of~$\W$,
we show that its moments of every fixed order $p\in\mathbb N$ are comparable to those of the
eigenvalue-trace statistic $Y_{k+1,n}=(k+1)^{-1}|\tr(\W^{k+1})|^2$. The resulting joint $(n,k)$-asymptotics reveal distinct scales for mean amplification and relative fluctuations of $Z_{k,n}$. The mean
$\E[Z_{k,n}]$ first departs from the classical fixed-$k$, large-$n$ prediction when the power scales as $k\asymp\sqrt n$. On the other hand, for each fixed integer $p\geq2$, the excess moment $\E[Z^p_{k,n}]/\E^p[Z_{k,n}] -1 $ becomes order one only when $k^2=\frac n p (\log n-\log\log n+O(1))$. In particular, when $\sqrt n\ll k\ll\sqrt{n\log n}$,
the mean diverges while relative fluctuations vanish: $Z_{k,n}/\mathbb E Z_{k,n}\to1$ in $L^p$ for every fixed $p\geq2$.}


\maketitle
\clearpage

\section{Introduction} \label{sec:intro}
Powers of nonnormal matrices\footnote{A complex square matrix $A$ is called nonnormal if $AA^*\neq A^*A$. } arise in many recurrent systems, including linear dynamical systems~\citep{trefethen2005spectra}, recurrent neural networks~\citep{orvieto2023resurrecting}, and state-space models~\citep{gu2021efficiently,gu2023mamba}. Consider the linear recurrence
\begin{equation}\label{eq:recurrence-intro}
    \mathbf h_{t+1}
    =
    \W \mathbf h_t + \mathbf B_n \mathbf u_t,
\end{equation}
where $\mathbf h_t$ is the recurrent state, $\mathbf u_t$ is the input, and
$\W$ is the state-transition matrix. The $k$-step Jacobian of this system is
\begin{equation}\label{eq:recurrent-jacobian-intro}
    J_{t,k}
    :=
    D_{\mathbf h_t}\mathbf h_{t+k}
    =
    \W^k.
\end{equation}
Thus powers of the matrix $\W$ govern the propagation of past inputs, effect of perturbations, and sensitivity to noise. Our main object of interest is the average squared singular value of the $k$-step Jacobian,
\begin{equation}\label{eq:Z-intro}
    Z_{k,n}
    :=
    \frac1n\|J_{t,k}\|_F^2
    =
    \frac1n\sum_{j=1}^n s_j(\W^k)^2 .
\end{equation}
Equivalently, $Z_{k,n}$ is the mean squared gain of $J_{t,k}$ over isotropic directions, and therefore gives a natural measure of average signal amplification over a lag of $k$ steps. In this work, we consider a normalized complex
Ginibre state-transition matrix:
\begin{equation}\label{eq:ginibre}
    \W=\frac{1}{\sqrt n}G_n,
    \qquad
    (G_n)_{ij}\overset{\mathrm{i.i.d.}}{\sim}\mathcal N_{\mathbb C}(0,1).
\end{equation}
This setting is motivated by recurrent neural networks at random Gaussian initialization, known as Glorot initialization~\citep{glorot2010understanding}. The $n^{-1/2}$ factor in $\W$ in this initialization aims to preserve the expected squared signal magnitude at a single step at large $n$. Precisely, it ensures that $\lim_{n\to\infty}\mathbb{E}[\|\W \mathbf u\|^2]= \|\mathbf u\|^2$. More generally, the singular values of $J_{t,k}$ are a central object of study in neural networks theory, since they quantify the amplification or attenuation of signals and backpropagated gradients over $k$ steps~\citep{schoenholz2017deep,chen2018dynamical,yang2019wide,seleznova:2026}.  

\subsection{\texorpdfstring{Double-Scaling Limit $n\to\infty$, $k=k(n)\to\infty$}{Double-Scaling Limit}}

In this work, we investigate the behavior of $Z_{k,n}$ in the double-scaling limit, where the matrix dimension $n$ and the power $k=k(n)$ tend to infinity simultaneously. This setting connects two classical asymptotic regimes: (i) the long-time limit $k\to\infty$ at fixed dimension $n$, commonly studied in dynamical systems, and (ii) the large-dimensional limit $n\to\infty$ at fixed power $k$, commonly studied in random matrix theory. As we explain below, these two regimes exhibit strikingly different behavior.

\noindent\textbf{Large $k$, fixed $n$ behavior.} At fixed dimension $n$, Gelfand's formula implies that the asymptotic growth rate of $Z_{k,n}$ is determined by the spectral radius $\rho(\W)$ of $\W$: 
\begin{equation}
    Z_{k,n}^{1/(2k)} \longrightarrow \rho(\W),
    \qquad k\to\infty,
\end{equation}
holds for every realization.
Thus, the long-time behavior depends on the realization-dependent spectral radius and remains random. Moreover, the expectation of $Z_{k,n}$ eventually diverges:
\begin{equation}
    \mathbb E[Z_{k,n}]
    \geq
    \frac{1}{n^2}\mathbb E[\rho(\W)^{2k}]
    \longrightarrow\infty,
    \qquad k\to\infty.
\end{equation}
Indeed, for every fixed $n$, the spectral radius of a normalized Ginibre matrix $\W$ has unbounded support~\citep{rider2003limit,rider2014extremal}. Therefore, arbitrarily large realizations contribute to the growth of the expectation at sufficiently long time horizons.

\noindent\textbf{Large $n$, fixed $k$ behavior.}
In contrast, when the power $k$ is fixed and the dimension
$n\to\infty$, we have
\begin{equation}
    \mathbb E[Z_{k,n}]\longrightarrow 1,
    \qquad
    \operatorname{Var}(Z_{k,n})=O(n^{-2})\longrightarrow 0.
\end{equation}
These limits follow from known results on finite words of random matrices~\citep{alexeev:2010,dubach:2021}. Hence, for every fixed $k$, $Z_{k,n}$ converges in probability to the deterministic value $1$. This large-dimensional behavior is closely related to the infinite-width limit widely used in neural network theory, where random network observables often become deterministic and analytically tractable~\citep{yang2019wide}.

\noindent\textbf{Jointly growing $k$ and $n$.}
The two classical limits therefore lead to qualitatively different conclusions. At fixed dimension, the long-time behavior of $Z_{k,n}$ is realization-dependent and its expectation diverges. At fixed power, the large-dimensional limit instead yields a deterministic value of one. Neither limit, however, describes the behavior of a large recurrent system over a time horizon that grows with its dimension.

Recent works have begun to investigate this double-scaling setting. \citet{bar2026revisiting} established an initial lower bound on $\mathbb E[Z_{k,n}]$, while \citet{seleznova:2026} derived an exact finite-$(n,k)$ expression for the expectation and characterized its joint scaling regimes. In particular, the latter work showed that the expectation departs from its infinite-$n$, fixed~$k$ limit at the critical scale $k\asymp\sqrt n$. However, these results are limited to the expectation and leave open the behavior of higher moments and fluctuations. In particular, it remains unclear whether the growth of $\mathbb E[Z_{k,n}]$ reflects typical amplification or is driven by rare large realizations. This motivates the central question of our work:
\begin{center}
\emph{How does the distribution of $Z_{k,n}$ behave for growing $k=k(n)$ as $n\to\infty$?}
\end{center}
More specifically, we seek to identify the joint scaling regimes in which higher moments depart from their infinite-$n$ limits, fluctuations become non-negligible, and rare amplification events begin to influence the behavior of $Z_{k,n}$. 

\noindent\textbf{Products of independent random matrices.} 
Related double-scaling limits have been studied extensively for products of independent random matrices, where both the number of factors $m$ and their dimension $n$ tend to infinity. In particular, \citet{liu2023lyapunov,liu2024phase,akemann2020universality} identified phase transitions in the local singular-value and eigenvalue statistics of independent Ginibre matrix products, governed by the ratio $m/n$. Thus, the transition from the fixed-$m$, infinite-$n$ limit occurs at $m\asymp n$ for the spectral statistics of independent matrix products, in contrast to the $k\asymp\sqrt n$ scale for the expected amplification of Ginibre matrix powers.
Joint $(m,n)$-limits have also been investigated in the context of deep feedforward neural networks, where Jacobians involve products of independent random weight matrices, possibly interleaved with nonlinear activations. These works have characterized fluctuations of propagated signals and gradients, as well as finite-$n$ corrections to infinite-$n$ behavior~\citep{hanin2020products,hanin2018neural,hanin2019finite,seleznova2022neural}. 

\subsection{Singular-Eigenvalue Bridge}

A central difficulty in understanding the behavior of $Z_{k,n}$ stems from the nonnormality of the Ginibre matrix $\W$. For normal matrices, the Frobenius norm of a matrix power is determined entirely by its eigenvalues. For nonnormal matrices, however, the eigenvalue spectrum alone may provide little information about the norms of matrix powers~\citep{trefethen1997pseudospectra}. Indeed, Ginibre eigenvectors are highly nonorthogonal~\citep{chalker1998eigenvector,bourgade2020distribution}, which can lead to substantial differences between eigenvalue- and singular-value statistics. This is apparent even without taking powers: for normalized Ginibre matrix $\W$, the spectral radius converges to $1$~\citep{rider2003limit,rider2014extremal}, whereas the largest singular value converges to $2$~\citep{bai:silverstein:2010}.

While Ginibre eigenvalue statistics are exceptionally tractable and admit explicit finite-$n$ correlation functions~\citep{ginibre1965statistical,byun2025progress,dubach:2018}, the singular-value statistics of matrix powers are considerably harder to analyze. Motivated by this contrast, we connect the singular-value statistic $Z_{k,n}$ to a more tractable statistic of the eigenvalues. Specifically, we introduce
\begin{equation} \label{eq:def_Y}
    Y_{M,n}
    :=
    \frac{1}{M}\left|\tr(\W^M)\right|^2
    =
    \frac{1}{M}\sum_{i,j=1}^n
    \lambda_i^M\overline{\lambda_j^M}, \quad M\in\mathbb N,
\end{equation}
where $\lambda_1,\ldots,\lambda_n$ are the eigenvalues of $\W$. Unlike $Z_{k,n}$, the moments of $Y_{M,n}$ can be expressed using the explicit Ginibre eigenvalue correlation functions. Moreover, \citet{dubach:2024} derived combinatorial formulas for its first two moments and connected them to the distribution of certain random permutation cycles.

The key insight of our work is that, despite their different spectral nature, the moments of $Z_{k,n}$ and $Y_{k+1,n}$ admit closely related combinatorial representations. By analyzing their Wick expansions, we establish the exact identity
\begin{equation}
    \mathbb E[Z_{k,n}]
    =
    \mathbb E[Y_{k+1,n}],
\end{equation}
and, more generally, the following bounds for every integer $p\geq2$:
\begin{equation}
    \label{eq:moment-bridge-intro}
    \frac{1}{p^p}\mathbb E[Y_{k+1,n}^p]
    \leq
    \mathbb E[Z_{k,n}^p]
    \leq
    \frac{1}{p!}\mathbb E[Y_{k+1,n}^p].
\end{equation}
Remarkably, these bounds hold for arbitrary finite $n$ and $k$. Thus, although $Z_{k,n}$ and $Y_{k+1,n}$ need not be close for individual matrix realizations, all their moments are comparable.

Furthermore, the connection extends to the fluctuations of $Z_{k,n}$ around its mean. To quantify these fluctuations, we introduce the \emph{excess moments}
\begin{equation}
    \label{eq:excess-moments-intro}
    R^{(p)}_{k,n}
    :=
    \frac{\mathbb E[Z_{k,n}^p]}{(\mathbb E Z_{k,n})^p}-1,
    \qquad
    q^{(p)}_{M,n}
    :=
    \frac{\mathbb E[Y_{M,n}^p]}{(\mathbb E Y_{M,n})^p}-p!, \quad k,M\in\mathbb N.
\end{equation}
Here, $1$ and $p!$ are the respective limiting normalized moments in the fixed-$k$, large-$n$ regime, where $Z_{k,n}$ concentrates around its mean while $Y_{k+1,n}/\mathbb E[Y_{k+1,n}]$ converges to a unit-mean exponential random variable. Building on the same combinatorial representations, we establish the following bounds
\begin{equation}
    \label{eq:excess-bridge-intro}
    \frac{1}{p^p}q^{(p)}_{k+1,n}
    \leq
    R^{(p)}_{k,n}
    \leq
    \frac{1}{p!}q^{(p)}_{k+1,n}.
\end{equation}
Together, the inequalities~\eqref{eq:moment-bridge-intro} and \eqref{eq:excess-bridge-intro} constitute our \emph{singular--eigenvalue bridge}: they allow us to determine the asymptotic growth of the moments of $Z_{k,n}$ and the scaling of its relative fluctuations, using the more tractable eigenvalue statistic $Y_{k+1,n}$.

\subsection{A Hierarchy of Double Scaling Regimes}

Although the singular-eigenvalue bridge does not determine the moments of $Z_{k,n}$ exactly, it is sufficiently sharp to identify their joint $(n,k)$-scaling regimes. Indeed, for each fixed $p$, the comparison constants above depend only on $p$, and not on $n$ or $k$. Therefore, the moments and excess moments of $Z_{k,n}$ exhibit the same asymptotic growth scales as their counterparts for $Y_{k+1,n}$. This allows us to establish the following hierarchy of joint $(n,k)$-scaling regimes.

For the raw moments of $Z_{k,n}$, we identify three distinct regimes. For every fixed integer $p\geq1$,
\begin{equation}
\label{eq:raw-moment-regimes}
\mathbb E[Z_{k,n}^p]
\longrightarrow
 \begin{cases}
1,
& k/\sqrt n\longrightarrow0,\\[2mm]
\displaystyle
\left(\frac{2\sinh(c^2/2)}{c^2}\right)^p,
& k/\sqrt n\longrightarrow c\in(0,\infty),\\[3mm]
+\infty,
& k/\sqrt n\longrightarrow+\infty.
\end{cases}
\end{equation}
The limits for $\mathbb E[Z_{k,n}^p]$ are established in Corollary~\ref{corr:critical_convergence_Zp}. The same three limits hold for $\mathbb E[Y_{k+1,n}^p]/p!$, as follows from Proposition~\ref{prop:bridge} and Theorem~\ref{thm:q-asymptotics}. Thus, the fixed-$k$, large-$n$ behavior persists as long as $k=o(\sqrt n)$, while a nontrivial limiting amplification emerges at $k\asymp\sqrt n$. We note that the first moment case $p=1$ was previously obtained by \citet{seleznova:2026}.

Remarkably, the excess moments $R^{(p)}_{k,n}$ exhibit a different hierarchy of joint scaling regimes. To describe it, define $s_{p,n}$ for every fixed integer $p\geq1$ by
\begin{equation}\label{eq:fluc_scaling}
\frac{k^2}{n} = \frac{1}{p}\log n - \frac{1}{p}\log\log n + s_{p,n}.
\end{equation}
For $p\geq2$, the behavior of the excess moments is given by
\begin{equation}
\label{eq:excess-moment-regimes}
\begin{cases}
R^{(p)}_{k,n}\longrightarrow0,
& s_{p,n}\longrightarrow-\infty,\\[2mm]
R^{(p)}_{k,n}=\Theta_{p,s}(1),
& s_{p,n}\longrightarrow s\in\mathbb R,\\[2mm]
R^{(p)}_{k,n}\longrightarrow+\infty,
& s_{p,n}\longrightarrow+\infty.
\end{cases}
\end{equation}
In particular, the critical scale for the $p$-th excess moment is $k\asymp\sqrt{n\log n}$, substantially larger than the $k\asymp\sqrt n$ scale governing the raw moments. Consequently, there is an intermediate regime
\begin{equation}
\sqrt n\ll k\ll\sqrt{n\log n},
\end{equation}
in which the expectation of $Z_{k,n}$ diverges while the relative variance vanishes. In this regime, $Z_{k,n}$ remains concentrated around its mean, so the growth of its moments reflects typical, self-averaging amplification rather than a breakdown of concentration.

\begin{figure}[t]
\centering
\includegraphics[width=1.00\linewidth]{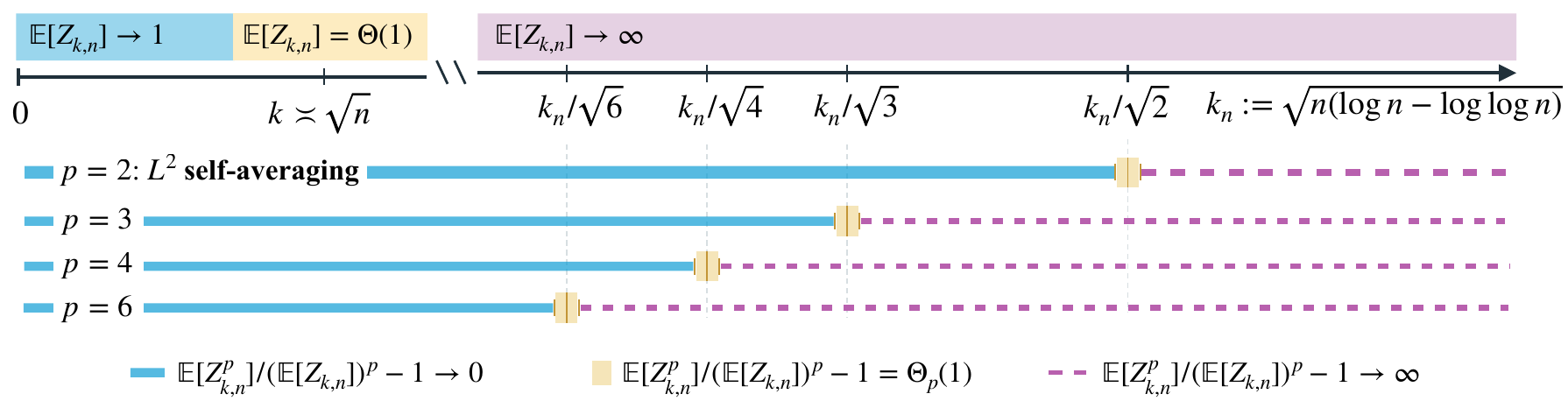}
\caption{Joint $(n,k)$-scaling regimes for the mean $\mathbb E[Z_{k,n}]$ (top band) and the $p$-th relative excess moments $R^{(p)}_{k,n}$ (lower bands, one for each $p$).}
\label{fig:scaling-hierarchy}
\end{figure}

Note that, by~\eqref{eq:fluc_scaling}, the excess moment critical scaling thresholds also depend on the moment order~$p$. Since the transition occurs when $k^2/n\sim(\log n)/p$, higher-order excess moments become significant at smaller values of $k$ than lower-order ones. This hierarchy reveals \emph{intermittency}: higher-order fluctuations become sensitive to increasingly rare, large amplification events even while lower-order fluctuations remain negligible. Figure~\ref{fig:scaling-hierarchy} summarizes the resulting joint $(n,k)$-scaling regimes.

\subsection{Outline and Notation}

The remainder of the paper is organized as follows. Section~\ref{sec:2} presents the main results and their numerical illustrations. Section~\ref{sec:proof-ideas} gives an overview of the proof structure and main ideas. Section~\ref{sec:rnns} applies the results to signal propagation in a linear recurrent network, and Section~\ref{sec:discussion} discusses future directions. The detailed proofs follow: Section~\ref{sec:proof-bridge} establishes the mean formula, its asymptotics, and the singular--eigenvalue bridge; Section~\ref{sec:proof-q-asymptotics} establishes the scaling regimes; and Section~\ref{sec:proof_thm_radial} proves the radial-dominance theorem used in the scaling analysis.

\vspace{1ex}
\noindent\textbf{Notation.}
We write $\mathcal{N}_{\mathbb{C}}(0,\sigma^2)$ for the proper complex Gaussian distribution, i.e.,  $Z\sim\mathcal{N}_{\mathbb{C}}(0,\sigma^2)$ means that
$Z=(\sigma/\sqrt{2})(X+\mathrm{i}Y)$,
where $X,Y$ are independent standard real Gaussian random variables.
For a matrix $\mathbf{A}\in\mathbb{C}^{n\times n}$, we denote its Frobenius norm by $\|\mathbf{A}\|_F$.
For positive sequences $a_n$ and $b_n$, we write $a_n\sim b_n$ if $a_n/b_n\to1$, $a_n\ll b_n$ if $a_n/b_n\to0$, and $a_n=\Theta(b_n)$ or $a_n\asymp b_n$ if there exist constants $0<c<C<\infty$ such that $cb_n\leq a_n\leq Cb_n$ for all sufficiently large $n$.
Unless stated otherwise, all asymptotic statements refer to $n\to\infty$, with $k=k(n)$ allowed to depend on $n$, while the moment order $p$ is fixed.
When comparing $Z_{k,n}$ and $Y_{M,n}$, we write $M=k+1$. Unless stated otherwise, the constants in asymptotic bounds may depend on $p$, but not on $n$ or $k$. %

\section{Main Results}\label{sec:2}

In this section, we state our main results formally. Throughout, $\mathbf W\in\mathbb C^{n\times n}$ is a normalized complex Ginibre matrix with independent entries distributed as $\mathcal N_{\mathbb C}(0,1/n)$. Recall the statistics $Z_{k,n}=n^{-1}\|\mathbf W^k\|_F^2$ and $Y_{M,n}=M^{-1}|\operatorname{tr}(\mathbf W^M)|^2$ introduced in the previous section. We establish the singular–eigenvalue bridge between their moments and excess moments, and characterize the resulting joint scaling regimes. Table~\ref{tab:regimes} summarizes the scaling regimes of the raw moments and excess moments of $Z_{k,n}$.

\begin{table}[!ht]
    \centering
    \begin{tabular}{c|c|c|c}
        \hline
        Regime
        & $\mathbb{E}[Z_{k,n}]$
        & $\mathbb{E}[Z_{k,n}^p]$
        & $R_{k,n}^{(p)}$
        \\
        \hline
        $\displaystyle \frac{k}{\sqrt n}\to0$
        & $\displaystyle \to1$
        & $\displaystyle \to1$
        & $\displaystyle \to0$
        \\[1.2ex]

        $\displaystyle \frac{k}{\sqrt n}\to c\in(0,\infty)$
        & $\displaystyle
          \to\frac{2}{c^2}
          \sinh\left(\frac{c^2}{2}\right)$
        & $\displaystyle
          \to\left(
          \frac{2}{c^2}
          \sinh\left(\frac{c^2}{2}\right)
          \right)^p$
        & $\displaystyle \to0$
        \\[2ex]

        $\displaystyle
        \frac{k}{\sqrt n}\to\infty,
        \quad s_{p,n}\to-\infty$
        & $\displaystyle
          \sim
          \frac{n}{k^2}
          \exp\left(\frac{k^2}{2n}\right)$
        & $\displaystyle
          \sim
          \left(
          \frac{n}{k^2}
          \exp\left(\frac{k^2}{2n}\right)
          \right)^p$
        & $\displaystyle \to0$
        \\[2ex]

        $\displaystyle s_{p,n}\to s\in\mathbb{R}$
        & $\displaystyle
          \sim
          \frac{n}{k^2}
          \exp\left(\frac{k^2}{2n}\right)$
        & $\displaystyle
          \Theta\left(
          \left(
          \frac{n}{k^2}
          \exp\left(\frac{k^2}{2n}\right)
          \right)^p
          \right)$
        & $\displaystyle \Theta(1)$
        \\[2ex]

        $\displaystyle s_{p,n}\to+\infty$
        & $\displaystyle \to+\infty$
        & $\displaystyle \to+\infty$
        & $\displaystyle \to+\infty$
        \\
        \hline
    \end{tabular}
    \caption{Joint $(n,k)$-scaling regimes for the mean
    $\mathbb{E}[Z_{k,n}]$, the $p$-th raw moment
    $\mathbb{E}[Z_{k,n}^p]$, and the $p$-th excess moment
    $R_{k,n}^{(p)}$, for fixed $p\geq2$. Here $s_{p,n}$ is defined in~\eqref{eq:fluc_scaling}.
    The mean asymptotics are given in
    Proposition~\ref{prop:bridge} and follow from
    \citet{seleznova:2026}. The excess-moment regimes
    are established in Theorem~\ref{thm:q-asymptotics},
    and the raw moment asymptotics follow from
    Corollary~\ref{corr:critical_convergence_Zp}.}
    \label{tab:regimes}
\end{table}

\subsection{Singular--Eigenvalue Bridge}

We first recall the mean behavior, which provides a starting point for the higher-moment analysis. The following proposition combines the identity $\mathbb E[Z_{k,n}]=\mathbb E[Y_{k+1,n}]$ with the mean asymptotics established by \citet{seleznova:2026}. We include the exact expression and asymptotics for completeness and for use in the subsequent scaling proofs.

\begin{proposition}\label{prop:bridge}
For all integers
$n,k\geq1$,
\begin{equation} \label{eq:mean-identity}
    \mathbb{E}[Z_{k,n}]
    =
    \mathbb{E}[Y_{k+1,n}].
\end{equation}
For $M=k+1$, their common value is given by
\begin{align}
    \mathbb E[Y_{M,n}]
    =
    \frac{n}{M(M+1)}
    \left[
        \prod_{j=0}^{M}\left(1+\frac jn\right)
        -
        \prod_{j=0}^{M}\left(1-\frac jn\right)
    \right].
    \label{eq:mean_expression}
\end{align}
Moreover, as $n\to\infty$ with $k=k(n)$,
\begin{equation}
\label{eq:mean-regimes}
\mathbb E[Z_{k,n}]
=
\mathbb E[Y_{k+1,n}]
\longrightarrow
\begin{cases}
1,
& k/\sqrt n\to0,\\[2mm]
\displaystyle
\frac{2\sinh(c^2/2)}{c^2},
& k/\sqrt n\to c\in(0,\infty),\\[3mm]
+\infty,
& k/\sqrt n\to+\infty.
\end{cases}
\end{equation}
In particular, in the regime $\sqrt n\ll k\ll n^{2/3}$,
\begin{equation}\label{eq:mean-regimes2}
    \mathbb E[Z_{k,n}]
    =
    \mathbb E[Y_{k+1,n}]
    \sim
    \frac{n}{k^2}
    \exp\left(\frac{k^2}{2n}\right).
\end{equation}
\end{proposition}
The mean identity~\eqref{eq:mean-identity} is proved in Section~\ref{subsec:mean-identity-ideas}, while the exact formula and asymptotics are proved in Section~\ref{subsec:proof-mean}.

 Our main result extends this connection to higher moments and excess moments. Recall the excess moments $R_{k,n}^{(p)}$ and $q_{M,n}^{(p)}$ defined in~\eqref{eq:excess-moments-intro}.
In particular, for $p=2$, $R^{(2)}_{k,n}$ is the relative variance of $Z_{k,n}$.
The following theorem establishes the singular--eigenvalue bridge inequalities for both raw moments and excess moments.

\begin{theorem}\label{thm:bridge-variance}
For all integers $n,k\geq1$ and $p\geq2$,
\begin{align}
    \frac{1}{p^p}\mathbb{E}[Y_{k+1,n}^p]
    &\leq
    \mathbb{E}[Z_{k,n}^p]
    \leq
    \frac{1}{p!}\mathbb{E}[Y_{k+1,n}^p],
    \label{eq:thm:bridge-variance1}
    \\
    \frac{1}{p^p}q^{(p)}_{k+1,n}
    &\leq
    R^{(p)}_{k,n}
    \leq
    \frac{1}{p!}q^{(p)}_{k+1,n}.
    \label{eq:thm:bridge-variance2}
\end{align}
\end{theorem}

The proof idea is explained in Section~\ref{sec:proof_ideas_bridge}, and the complete proof is given in Section~\ref{sec:proof-bridge}.

\subsection{Joint Scaling Regimes}

The comparison constants in Theorem~\ref{thm:bridge-variance} depend only on $p$, and not on $n$ or $k$. Therefore, we can determine the joint scaling regimes of $R^{(p)}_{k,n}$ by studying $q^{(p)}_{k+1,n}$. 
The following theorem describes the scaling regimes of the excess moments
in terms of $s_{p,n}$ defined in~\eqref{eq:fluc_scaling}.

\begin{theorem}\label{thm:q-asymptotics}
Fix an integer $p\geq2$. As $n\to\infty$ with
$k=k(n)$, the following regimes hold:
\begin{enumerate}
    \item
    If $s_{p,n}\to-\infty$, then
    \[
        q^{(p)}_{k+1,n}\longrightarrow0,
        \qquad
        R^{(p)}_{k,n}\longrightarrow 0.
    \]
    In particular,
    \[
        \frac{Z_{k,n}}{\mathbb{E}Z_{k,n}}
        \longrightarrow1
        \qquad\text{in }L^p.
    \]

    \item 
    If $s_{p,n}\to s\in\mathbb{R}$, then
    \[
        q^{(p)}_{k+1,n}
        \longrightarrow
        p^{-(p+1)/2}
        \exp\left(\frac{p(p-1)}{2}s\right).
    \]
    Consequently,
    \[
        R^{(p)}_{k,n}=\Theta(1).
    \]

    \item
    If $s_{p,n}\to+\infty$, then
    \[
        q^{(p)}_{k+1,n}\longrightarrow+\infty,
        \qquad
        R^{(p)}_{k,n}\longrightarrow+\infty.
    \]
\end{enumerate}
\end{theorem}
In particular, if $k/\sqrt n\to\infty$ and $s_{2,n}\to-\infty$, then
\[
 \mathbb E Z_{k,n}\longrightarrow\infty,
 \qquad \frac{Z_{k,n}}{\mathbb E Z_{k,n}}\longrightarrow1
 \quad\text{in }L^2.
\]
The main ideas behind Theorem~\ref{thm:q-asymptotics} are summarized in Section~\ref{sec:proof_ideas_scaling}. The complete proof is given in Section~\ref{sec:proof-q-asymptotics}.

Combining these excess-moment limits with the mean asymptotics in Proposition~\ref{prop:bridge}, we obtain the following joint scaling regimes for the raw moments.

\begin{corollary}\label{corr:critical_convergence_Zp}
Fix an integer $p\geq1$. As $n\to\infty$
with $k=k(n)$,
\begin{equation}
\label{eq:raw-moment-regimes-detail}
\mathbb E[Z_{k,n}^p]
\longrightarrow
\begin{cases}
1,
& k/\sqrt n\longrightarrow0,\\[2mm]
\displaystyle
\left(\frac{2\sinh(c^2/2)}{c^2}\right)^p,
& k/\sqrt n\longrightarrow c\in(0,\infty),\\[3mm]
+\infty,
& k/\sqrt n\longrightarrow+\infty.
\end{cases}
\end{equation}
Moreover, if $k\gg\sqrt{n}$ and $s_{p,n}\to -\infty$,
\begin{equation}
    \mathbb E[Z_{k,n}^p]
    \sim
    \left(
        \frac{n}{k^2}
        \exp\left(\frac{k^2}{2n}\right)
    \right)^p.
\end{equation}
\end{corollary}

The proof is given in Section~\ref{subsec:proof-raw}.

\subsection{Numerical Illustrations}\label{subsec:numerical-illustrations}

\begin{figure}[!ht]
    \centering
    \includegraphics[width=\linewidth]{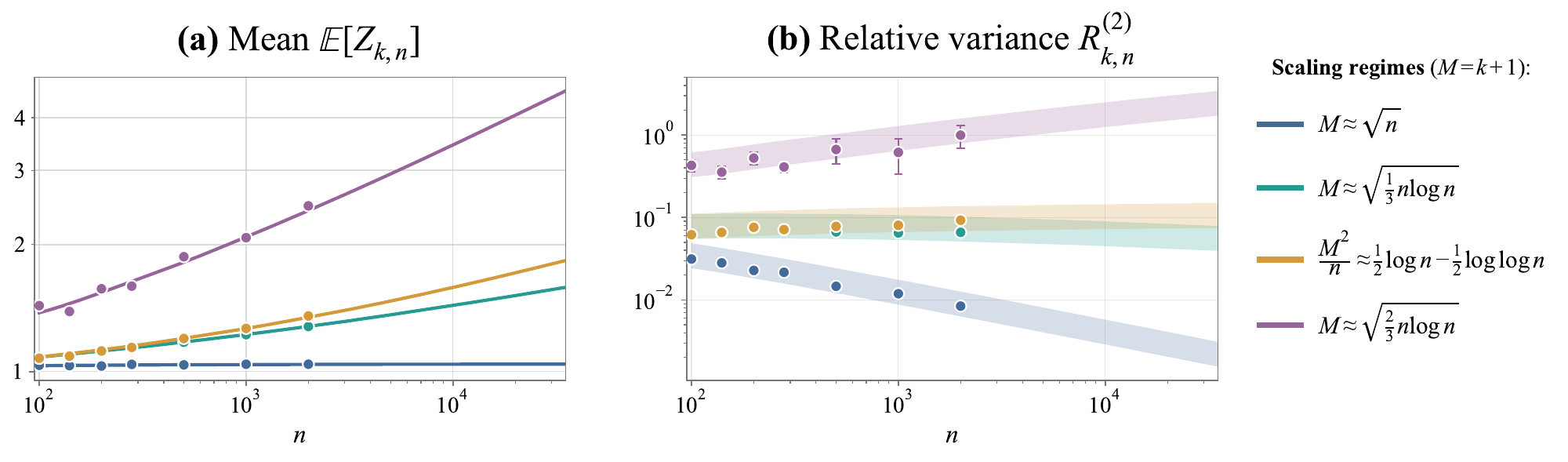}
    \caption{\textbf{Mean growth and relative fluctuations occur on different scales.} We compute $Z_{k,n}=n^{-1}\|\W^k\|_F^2$ for $n\in\{100,140,200,280,500,1000,2000\}$, with $M=k+1$ chosen according to the four scaling regimes specified in the legend. Panel \textbf{(a)} shows mean $\mathbb{E}[Z_{k,n}]$, panel \textbf{(b)} shows the relative variance $R^{(2)}_{k,n} =\operatorname{Var}(Z_{k,n})/(\mathbb{E}[Z_{k,n}])^2$. Points show Monte Carlo estimates of the two quantities, using $10^4$ samples, with error bars indicating $\pm2$ estimated standard errors. The relative variance is estimated by averaging $(Z_{k,n}/\mu_{k,n}-1)^2$, using the exact finite-$n$ mean $\mu_{k,n}=\mathbb{E}[Z_{k,n}]$. Solid curves in~\textbf{(a)} show the smoothed theoretical mean value. Shaded bands in~\textbf{(b)} show the smoothed bounds $q^{(2)}_{M,n}/4\leq R^{(2)}_{k,n}\leq q^{(2)}_{M,n}/2$. The intermediate regime $M\approx\sqrt{n\log n/3}$ illustrates that the mean can grow while relative fluctuations decrease.}
    \label{fig:mean_variance}
\end{figure}

We now illustrate the joint scaling regimes above with numerical experiments. Figure~\ref{fig:mean_variance} compares the mean and relative variance of $Z_{k,n}$ with their theoretical values and bounds. In particular, the intermediate regime shows that the mean can grow while the relative fluctuations decrease, as predicted by Proposition~\ref{prop:bridge} and Theorem~\ref{thm:q-asymptotics}.
\begin{figure}[h!]
    \centering
    \includegraphics[width=\linewidth]{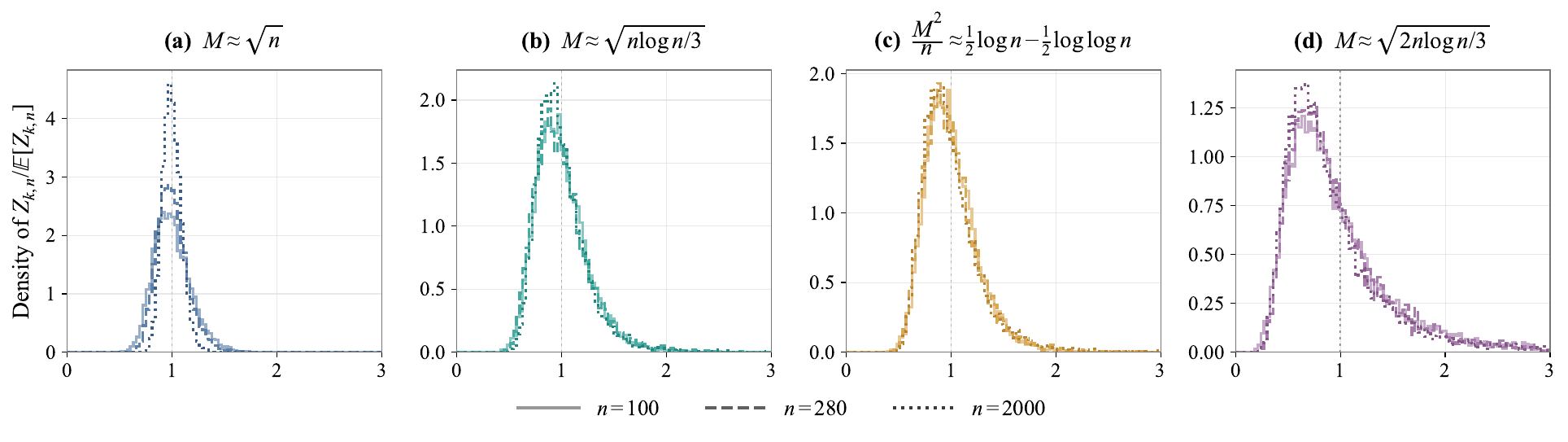}
    \caption{\textbf{Density of $Z_{k,n}/\mathbb E[Z_{k,n}]$ across scaling regimes.} Panels show empirical densities of $Z_{k,n}/\mu_{k,n}$, where $\mu_{k,n}=\mathbb{E}[Z_{k,n}]$ is the exact theoretical mean, in the four scaling regimes from Fig.~\ref{fig:mean_variance}, specified in the panel titles. Line styles and shades distinguish the dimensions $n$ listed in the legend. In~\textbf{(a)}, increasing $n$ leads to concentration around the normalized mean, marked by the vertical line at~$1$. Regime~\textbf{(b)} also exhibits asymptotic concentration, but convergence is much slower. In~\textbf{(c)} and~\textbf{(d)}, substantial right tails persist,
    and there is no concentration in the relative $L^2$ sense.}
    \label{fig:normalized-histograms}
\end{figure}

To further illustrate the behavior of the fluctuations, Figure~\ref{fig:normalized-histograms} shows empirical densities of $Z_{k,n}/\mathbb E Z_{k,n}$ in the same four regimes.
In panels (a) and (b), both of which are subcritical for the relative
variance, the empirical distributions become increasingly concentrated
around one as $n$ grows. This concentration is less pronounced in
panel (b), where the powers are larger. In panels (c) and (d), which
show the critical and supercritical regimes for the relative variance,
respectively, the empirical distributions exhibit substantial right tails.

\section{Overview of the Proofs}\label{sec:proof-ideas}

In this section, we give an overview of the proof structure and explain the main ideas behind the singular–eigenvalue bridge and the joint scaling regimes. The detailed arguments are given in Sections~\ref{sec:proof-bridge} and~\ref{sec:proof-q-asymptotics}.

\subsection{Singular-Eigenvalue Bridge}\label{sec:proof_ideas_bridge}
We start by summarizing the proof idea for the singular–eigenvalue bridge inequalities in Theorem~\ref{thm:bridge-variance}. Before considering general $p$, we discuss the special case $p=1$, stated in Proposition~\ref{prop:bridge}, which admits a particularly simple proof.

\subsubsection{\texorpdfstring{Mean Identity (The Case $p=1$)}{Mean Identity (The Case p=1)}}\label{subsec:mean-identity-ideas}
\begin{figure}[h]
    \centering
    \includegraphics[width=\linewidth]{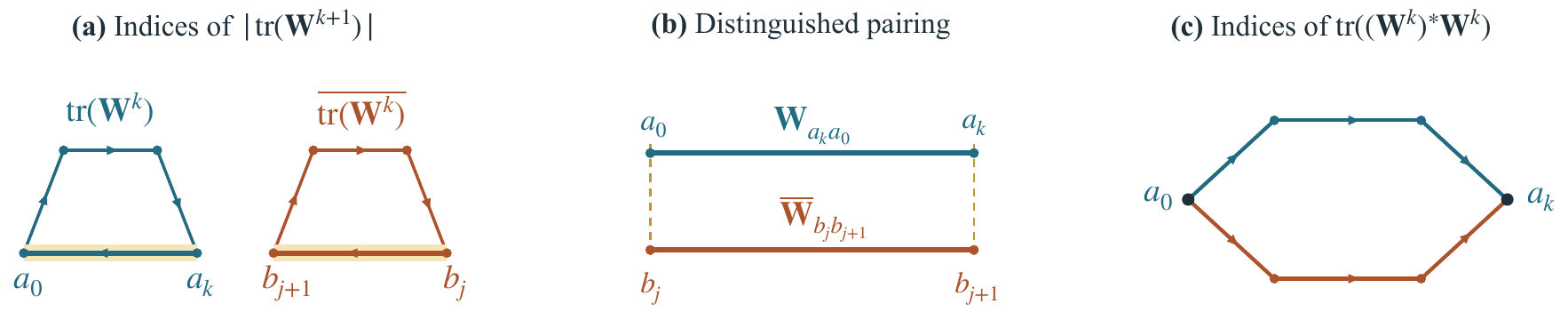}
    \caption{Illustration of the mean identity
    $\E|\tr(\W^{k+1})|^2
    =\frac{k+1}{n}\E\tr((\W^k)^*\W^k)$
    through a distinguished Wick pairing.
    Pairing one edge from each of the two cycles identifies
    their endpoints, leaving two paths of length $k$
    with common endpoints.}
    \label{fig:wick-bridge}
\end{figure}

We first explain the mean identity
$\E Y_{k+1,n}=\E Z_{k,n}$ through Wick's formula,
using the graphical interpretation in
Figure~\ref{fig:wick-bridge}.
Notice that the index expansion of $|\tr(\W^{k+1})|^2$
involves two cycles of length $k+1$, corresponding
to the non-conjugate and conjugate matrix entries, as illustrated in Figure~\ref{fig:wick-bridge}(a). Since the entries of $\W$ are proper complex Gaussian variables with $\E[\W_{ab}\W_{cd}]=0$, every nonzero Wick pairing connects an edge of the first cycle to an edge of the second.
On the other hand, the expansion of
$\tr((\W^k)^*\W^k)$ involves two paths of
length $k$ with common endpoints, as in Figure~\ref{fig:wick-bridge}(c).

We connect these two structures via a distinguished Wick pairing.
Fix one edge in the first cycle and pair it with
an edge in the second cycle, as illustrated in
Figure~\ref{fig:wick-bridge}(b).
Gluing these edges leaves precisely the two paths in the Frobenius-norm expansion.
Since the distinguished edge in the first cycle has
$k+1$ possible partners in the second cycle, and the distinguished pairing's covariance contributes a factor $1/n$, this correspondence yields the factor $(k+1)/n$.

More precisely, writing $b_{k+1}=b_0$ and using
Wick's formula, we obtain
\begin{equation}\label{eq:mean_brige_proof}
\begin{aligned}
\E\bigl|\tr(\W^{k+1})\bigr|^2
&=
\sum_{\substack{a_0,\ldots,a_k\\b_0,\ldots,b_k}}
\E\Big[
\underbrace{\W_{a_0a_1}\W_{a_1a_2}\cdots \W_{a_ka_0}}_{\text{first cycle}}
\underbrace{\overline{\W_{b_0b_1}\W_{b_1b_2}\cdots \W_{b_kb_0}}}_{\text{second cycle}}
\Big]\\
&=
\sum_{j=0}^{k}
\sum_{\substack{a_0,\ldots,a_k\\b_0,\ldots,b_k}}
\underbrace{\E\left[
\W_{a_ka_0}\overline{\W_{b_jb_{j+1}}}
\right]}_{\text{distinguished pairing}}
\E\left[
\W_{a_0a_1}\cdots \W_{a_{k-1}a_k}
\prod_{\substack{r=0\\r\ne j}}^k
\overline{\W_{b_rb_{r+1}}}
\right]\\
&=
\frac{k+1}{n}
\sum_{\substack{a_0,\ldots,a_k\\b_1,\ldots,b_{k-1}}}
\E\left[
\W_{a_0a_1}\cdots\W_{a_{k-1}a_k}
\overline{\W_{a_0b_1}\cdots\W_{b_{k-1}a_k}}
\right]\\
&=
\frac{k+1}{n}\E\tr\bigl((\W^k)^*\W^k\bigr).
\end{aligned}
\end{equation}
Here the third equality follows by cyclic symmetry,
which allows us to take $j=k$, and by the covariance
identity
\[
\E\left[
\W_{a_ka_0}\overline{\W_{b_kb_0}}
\right]
=
\frac1n\delta_{a_kb_k}\delta_{a_0b_0}.
\]
The resulting identifications $b_k=a_k$ and
$b_0=a_0$ give the common endpoints of the two
remaining paths. Dividing \eqref{eq:mean_brige_proof}
by $k+1$ proves $\E Y_{k+1,n}=\E Z_{k,n}$.

The exact mean formula and asymptotics are proved in Section~\ref{subsec:proof-mean}.

\subsubsection{\texorpdfstring{Higher-Moment Bounds (The Case $p\geq2$)}{Higher-Moment Bounds (The Case p>=2)}}\label{sec:proof-ideas-higher-moments}

To establish the higher-moment bounds in Theorem~\ref{thm:bridge-variance}, we extend the
Wick pairing argument from the mean case to general
$p\geq2$. The key observation is that the moments of
$Y_{k+1,n}$ and $Z_{k,n}$ admit closely related
representations as weighted sums over permutations.

Set $M:=k+1$ and, for $i=1,\ldots,p$, introduce
the disjoint labelled sets $I_i:=\{0_i,\ldots,k_i\}$.
Let us fix the following permutation
\[
    \tau_p
    :=
    (0_1,1_1,\ldots,k_1)\cdots
    (0_p,1_p,\ldots,k_p)
    \in S\left(\bigcup_{i=1}^p I_i\right),
\]
and define the set of all permutations with
the same cycle structure as $\tau_p$:
\begin{equation}
    \mathcal C_{M,p}
    :=
    \left\{
        c\in S\left(\bigcup_{i=1}^p I_i\right):
        c \text{ consists of } p
        \text{ cycles, each of length } M
    \right\}.
\end{equation}
Here $S(A)$ denotes the permutation group of a finite
set $A$, and $\cyc(\pi)$ is the number of cycles of
$\pi$, including fixed points.

By Lemma~\ref{lem:y-perm}, which builds on the
permutation representation of Ginibre eigenvalue
statistics in \citet{dubach:2024}, we have
\begin{equation}
    \E[Y_{M,n}^p]
    =
    \frac{p!}{n^{pM}}
    \sum_{c\in\mathcal C_{M,p}}
    n^{\cyc(\tau_p^{-1}c)}.
\end{equation}
For the Frobenius-norm statistic, we obtain the
analogous representation
\begin{equation}
    \E[Z_{k,n}^p]
    =
    \frac{1}{n^{pM}}
    \sum_{c\in\mathcal C_{M,p}^{0}}
    n^{\cyc(\tau_p^{-1}c)},
\end{equation}
where
\begin{equation}
    \mathcal C_{M,p}^{0}
    :=
    \left\{
        c\in\mathcal C_{M,p}:
        0_1,\ldots,0_p
        \text{ belong to distinct cycles of }c
    \right\}.
\end{equation}
This representation is established in
Lemma~\ref{lem:Z_commutator_representation}.

Thus, both moments are expressed through the same
weights $n^{\cyc(\tau_p^{-1}c)}$, but the sums range
over different sets of permutations. The bridge
inequalities therefore reduce to comparing the
weighted sums over $\mathcal C_{M,p}$ and its
subset $\mathcal C_{M,p}^{0}$. This provides a generalization of the distinguished Wick pairing argument used for the mean. The representations are proved in Section~\ref{subsec:bridge-representations}.

The central comparison, established in
Lemma~\ref{lem:cycle_sum_comparison}, is
\begin{equation}\label{eq:C_sum_inequalities}
    \frac{p!}{p^p}
    \sum_{c\in\mathcal C_{M,p}}
        n^{\cyc(\tau_p^{-1}c)}
    \leq
    \sum_{c\in\mathcal C_{M,p}^{0}}
        n^{\cyc(\tau_p^{-1}c)}
    \leq
    \sum_{c\in\mathcal C_{M,p}}
        n^{\cyc(\tau_p^{-1}c)}.
\end{equation}
The main ingredient of the proof is the following.
For each $c\in\mathcal C_{M,p}$, consider the
proportion $w(c)$ of the $M^p$ choices $(u_1,\ldots,u_p)\in I_1\times\cdots\times I_p$
for which $u_1,\ldots,u_p$ belong to distinct
cycles of $c$. We show that this proportion satisfies
the uniform lower bound
\[
    w(c)\geq\frac{p!}{p^p}.
\]
The comparison is proved in Section~\ref{subsec:bridge-comparison}. This estimate, together with the permutation
representations above, yields the moment inequalities
in Theorem~\ref{thm:bridge-variance}. To obtain the corresponding inequalities for the excess moments, Corollary~\ref{cor:connected_cycle_sum_comparison} establishes the same comparison as in~\eqref{eq:C_sum_inequalities} with the sums restricted to the connected permutation sets
$\mathcal C_{M,p}^{\mathrm{conn}}$ and $\mathcal C_{M,p}^{0,\mathrm{conn}}$.

\subsection{Joint Scaling Regimes}\label{sec:proof_ideas_scaling}
The singular--eigenvalue bridge allows us to study the
moments and excess moments of $Z_{k,n}$ in \eqref{eq:Z-intro} through those of $Y_{M,n}$ in \eqref{eq:def_Y}.
The advantage of this reduction is that the joint
density of the Ginibre eigenvalues is known explicitly~\cite{ginibre1965statistical}.
Therefore, the moments of $Y_{M,n}$ can, in
principle, be computed by integration against this
density, making the joint $(M,n)$-scaling problem tractable.

However, for general moment order $p$,
the resulting integrals involve interactions among
multiple eigenvalues and are difficult to evaluate
in a form suitable for asymptotic analysis.
To overcome this difficulty, we show that the
excess moment asymptotics are largely governed by
a simpler radial term $A_{pM}(n)$, defined below,
which can be evaluated explicitly by integration
against the Ginibre eigenvalue density.

\subsubsection{Excess Moments via Mixed Cumulants}

We begin by expressing the excess moment $q^{(p)}_{M,n}$ of $Y_{M,n}$ via cumulants. Let $T:=\operatorname{tr}(\mathbf W^M)=\sum_j\lambda_j^M$ and define the balanced mixed cumulants
\[
    K_q=\kappa(
        \underbrace{T,\ldots,T}_{q},
        \underbrace{\overline T,\ldots,\overline T}_{q}),
    \qquad K_1=M\mathbb E[Y_{M,n}].
\]
Note that rotational invariance makes all mixed cumulants with unequal numbers of $T$ and $\overline T$ vanish. Therefore, only balanced cumulants appear in the moment-cumulant formula~\eqref{eq:moment-cumulant} expressing $q^{(p)}_{M,n}$. We can write the formula schematically as follows:
\[
    q^{(p)}_{M,n}
    =\frac{K_p}{K_1^p}
    +\text{terms involving lower-order normalized cumulants},
\]
where the remaining terms are products of normalized
cumulants $K_q/K_1^q$ with $2\leq q<p$. The exact expansion is given in Section~\ref{subsec:q-auxiliary}. The core of our proof scaling proofs is that, around the critical fluctuation regime, the scaling of excess moments is governed by $K_p/K_1^p$, and $K_p$ itself grows as the radial term $A_{pM}(n)$.

\subsubsection{Radial Dominance}

We define the radial term $A_{qM}(n)$ as the term
in the moment expansion that depends only on the
one-point Ginibre eigenvalue density.
Due to its simple structure, it admits an explicit
formula, derived in Lemma~\ref{lem:radial-term}:
\[
    A_{qM}(n)
    :=\mathbb E\!\left[\sum_{j=1}^n|\lambda_j|^{2qM}\right] = \frac{(n+qM)!}{(qM+1)(n-1)!\,n^{qM}}.
\]
Theorem~\ref{thm:radial-dominance} then establishes
that the difference between the balanced mixed
cumulant $K_q$ and $A_{qM}(n)$ satisfies
\[
    |K_q-A_{qM}(n)|
    \leq C_q A_{qM}(n)
    \exp\!\left(
        -\frac{(q-1)M^2}{n+qM}
    \right).
\]
In particular, when $M \gg \sqrt{n}$, the exponential factor vanishes, and hence
\[
    K_q\sim A_{qM}(n).
\]
Thus, in this regime, the scaling behavior of
the mixed cumulants is determined by an explicit
radial moment term. This reduction is central
to our analysis of the joint scaling regimes.
The proof of the radial dominance Theorem~\ref{thm:radial-dominance} is given in Section~\ref{sec:proof_thm_radial}.

\subsubsection{Intermediate Powers and the Critical Windows}
The excess moments become order one in the range
$\sqrt n\ll M\ll n^{2/3}$, which we refer to as
the intermediate powers regime. Here, the radial term formula and the mean asymptotics
from Proposition~\ref{prop:bridge} give
\[
    K_q\sim\frac{n}{qM}
    \exp\!\left(\frac{q^2M^2}{2n}\right),
    \qquad
    K_1\sim\frac nM
    \exp\!\left(\frac{M^2}{2n}\right).
\]
Consequently,
\[
    \frac{K_q}{K_1^q}
    \sim\frac1q\left(\frac Mn\right)^{q-1}
    \exp\!\left(\frac{q(q-1)M^2}{2n}\right).
\]
From this expression, in Section~\ref{subsec:q-intermediate-powers}, we deduce that the highest-order normalized cumulant $K_p/K_1^p$ becomes of order one precisely when
\[
    \frac{M^2}{n}
    =\frac{\log n-\log\log n}{p}+O(1).
\]
Moreover, in this critical regime, all lower-order
normalized cumulants $K_q/K_1^q$, $2\leq q<p$,
tend to zero. Hence the excess moment is asymptotically determined by the highest-order term:
\[
    q^{(p)}_{M,n}\sim\frac{K_p}{K_1^p}.
\]
In particular, if $s_{p,n}\to s\in\mathbb R$, then
\[
    q^{(p)}_{M,n}\longrightarrow
    p^{-(p+1)/2}
    \exp\!\left(\frac{p(p-1)}2s\right).
\]
Below this critical window, all normalized cumulants
up to order $p$ vanish, while above it the
order-$p$ contribution diverges.
This yields the three joint scaling regimes
of Theorem~\ref{thm:q-asymptotics}.

\subsubsection{The Remaining Regimes}
To complete the proof of Theorem~\ref{thm:q-asymptotics}, we must also consider the scaling regimes outside the intermediate powers range $\sqrt n\ll M\ll n^{2/3}$. The regime $M=O(\sqrt n)$, in which the excess moments vanish, is treated in Section~\ref{subsec:q-small-powers}, while the regime $M\geq c n^{2/3}$ is addressed in Section~\ref{subsec:q-large-powers}. Both cases require shorter arguments than the intermediate powers regime and are presented for completeness.

\section{Signal Propagation in Linear Recurrent Networks}\label{sec:rnns}

As a first step toward understanding signal propagation in deep recurrent neural networks, we consider a simplified linear recurrent network. 
Writing $\mathbf x_n=\mathbf h_0$ for the input, the network's output after $k+1$ recurrent steps is given by
\begin{equation}\label{eq:network_state}
    \mathbf h_{k+1,n}
    =\W\mathbf h_{k,n}
    =\W^{k+1}\mathbf x_n.
\end{equation}
We impose the following assumption on the signal $\mathbf x_n.$
\begin{assumption}\label{ass:input}
The random vector $\mathbf{x}_n=(x_1,\ldots,x_n)^\top$ is independent of $\mathbf{W}$. The entries $x_1, \ldots, x_n$ are i.i.d. real-valued random variables satisfying 
\begin{align*}
\E x_1=0,\qquad \E x_1^2=1,\qquad \mu_4:=\E x_1^4<\infty.
\end{align*}
\end{assumption}
The normalized signal energy
$n^{-1}\|\mathbf h_{k,n}\|^2$ measures the average squared
activation per neuron. Since
$\mathbb{E}\|\mathbf h_{0,n}\|^2=n$,
values of $n^{-1}\mathbb{E}\|\mathbf h_{k,n}\|^2$
above one indicate amplification of the initial signal
on average, whereas values below one indicate weakening
of the signal. The second moment of the normalized signal
energy, or equivalently the fourth moment of the hidden-state norm,
\begin{align} \label{eq:fourth_mom}
    \frac{1}{n^2} \E  \|\mathbf h_{k,n}\|^4 
\end{align}
allows us to quantify its fluctuations with respect
to the randomness of the weights and inputs.
Combining the results of Section \ref{sec:2} with standard moment bounds for quadratic forms gives the following result for \eqref{eq:fourth_mom}.
\begin{proposition}\label{prop:moments}
Suppose Assumption~\ref{ass:input} holds. Let $C_x:=\max\{3,\mu_4\}$.
Then, it holds for every $k,n\geq1$
\begin{equation}
1+\frac14q_{k+1,n}^{(2)}\leq\frac{\E\|\mathbf h_{k,n}\|^4}{(\E\|\mathbf h_{k,n}\|^2)^2}
\leq C_x\left(1+\frac12q_{k+1,n}^{(2)}\right).
\end{equation}
\end{proposition}
\label{sec:proof-moments}
\begin{proof}
    Let $B=\operatorname{Re}((\mathbf W^k)^*\mathbf W^k)$.
    Since $\mathbf x_n$ is real, we have
$\|\mathbf h_{k,n}\|^2=\mathbf x_n^\top B\mathbf x_n$.
Using a standard formula for the variance of quadratic forms (see, e.g., (9.8.6) in \cite{bai:silverstein:2010}), we obtain 
\begin{align*}
     \Var(\|\mathbf h_{k,n}\|^2\mid\mathbf W)
    =(\mu_4-1)\sum_i B_{ii}^2+2\sum_{i\ne j}B_{ij}^2
    \leq (C_x-1)\tr(B^2)
    \leq (C_x-1)(\tr B)^2.
\end{align*}
Since $\E[\|\mathbf h_{k,n}\|^2\mid\mathbf W]=\tr B$,
adding the squared conditional mean gives
\begin{align*}
    (\tr B)^2 \leq \E[\|\mathbf h_{k,n}\|^4\mid\mathbf W] \leq 
C_x(\tr B)^2
\end{align*}
Taking expectations and using $\tr B=nZ_{k,n}$
gives
\begin{align*}
    1+R^{(2)}_{k,n}
    \leq
    \frac{\E\|\mathbf h_{k,n}\|^4}
         {(\E\|\mathbf h_{k,n}\|^2)^2}
    \leq C_x(1+R^{(2)}_{k,n}).
\end{align*}
The assertion follows from \eqref{eq:thm:bridge-variance2} with $p=2$.
\end{proof}
Proposition~\ref{prop:moments} allows us to transfer the
fluctuation regimes of $Z_{k,n}$ established in
Theorem~\ref{thm:q-asymptotics} to moment bounds for the hidden state.
\begin{corollary}\label{cor:regimes}
Suppose Assumption~\ref{ass:input} holds. Then, the following assertions hold:
\begin{enumerate}
    \item If $s_{2,n}\to-\infty$, then
    \begin{align*}
        \E\|\mathbf h_{k,n}\|^4
        &=\Theta\!\left(
            (\E\|\mathbf h_{k,n}\|^2)^2
        \right), \qquad 
        \frac{\Var(\|\mathbf h_{k,n}\|^2)}
             {(\E\|\mathbf h_{k,n}\|^2)^2}
        =O(1).
    \end{align*}

    \item If $s_{2,n}\to s\in\mathbb R$, then
    \begin{align*}
        \E\|\mathbf h_{k,n}\|^4
        &=\Theta\!\left(
            (\E\|\mathbf h_{k,n}\|^2)^2
        \right),\qquad 
        \frac{\Var(\|\mathbf h_{k,n}\|^2)}
             {(\E\|\mathbf h_{k,n}\|^2)^2}
        =\Theta(1).
    \end{align*}

    \item If $s_{2,n}\to+\infty$, then
    \begin{align*}
        \E\|\mathbf h_{k,n}\|^4
        &\to\infty,\qquad 
        \frac{\Var(\|\mathbf h_{k,n}\|^2)}
             {(\E\|\mathbf h_{k,n}\|^2)^2}
        \to\infty.
    \end{align*}
\end{enumerate}
\end{corollary}
Corollary~\ref{cor:regimes} shows that the relative variance
of the signal energy remains bounded when $s_{2,n}\to-\infty$,
is of constant order when $s_{2,n}\to s\in\mathbb{R}$, and
diverges when $s_{2,n}\to+\infty$.

Moreover, the upper bound in Proposition \ref{prop:moments}
gives a lower bound on the probability that the signal energy
exceeds a fixed fraction of its expectation: By the Paley-Zygmund inequality, we have for every $\theta\in(0,1)$,
\begin{align}
    \mathbb{P}\!\left(
        \|\mathbf h_{k,n}\|^2
        \geq \theta\,\mathbb{E}\|\mathbf h_{k,n}\|^2
    \right)
    &\geq
    (1-\theta)^2
    \frac{(\mathbb{E}\|\mathbf h_{k,n}\|^2)^2}
         {\mathbb{E}\|\mathbf h_{k,n}\|^4}
    \geq
    \frac{(1-\theta)^2}
         {C_x(1+q_{k+1,n}^{(2)}/2)}.
    \label{eq:signal_amplification_probability}
\end{align}
Consider the regime $k/\sqrt n\to\infty$ and
$s_{2,n}\to-\infty$. 
Theorem~\ref{thm:q-asymptotics} implies
$q_{k+1,n}^{(2)}\to0$, and
$n^{-1}\mathbb{E}\|\mathbf h_{k,n}\|^2\to\infty$ by Proposition~\ref{prop:bridge}.
Thus, the normalized signal energy
$n^{-1}\|\mathbf h_{k,n}\|^2$ exceeds any fixed
level with a probability bounded away from zero
as $n\to\infty$.
The growth of its expectation is therefore not caused
solely by large signal energies occurring with
vanishing probability.

\section{Discussion and Future Work}\label{sec:discussion}

Double-scaling limits have received considerable attention for products of independent random matrices~\citep{liu2023lyapunov,liu2024phase,akemann2020universality}, while the corresponding problem for powers of a single random matrix remains less understood. Building on previous results for the mean $\E Z_{k,n}$~\citep{bar2026revisiting,seleznova:2026}, this work provides, to the best of our knowledge, one of the first analyses of singular-value statistics in this setting. Our results reveal distinct scales for mean amplification and relative fluctuations of $Z_{k,n}$, together with a hierarchy of critical windows depending on the moment order.

The complex Ginibre ensemble provides a tractable starting point for this largely open problem. Extending the results to real Ginibre matrices and determining whether the scaling hierarchy persists for non-Gaussian entries are natural next questions. In this direction, \citet{seleznova:2026} showed that the mean of $Z_{k,n}$ for real Ginibre matrices is bounded below by its complex Ginibre counterpart, indicating that the scale $k\asymp\sqrt n$ already
leads to mean amplification in the real case.
Whether the higher-moment and fluctuation regimes
extend to this setting remains open.
More generally, establishing universality in double-scaling regimes requires estimates that remain
uniform as the power grows with the dimension.

For applications to recurrent systems, including neural networks, an important further direction is to incorporate nonlinearities between successive applications of the recurrent weight matrix. Understanding how these nonlinearities affect the amplification and fluctuation regimes would extend the analysis toward signal propagation and gradient dynamics in nonlinear recurrent systems, going beyond existing work on signal propagation
in deep recurrent networks~\citep{schoenholz2017deep,chen2018dynamical}.

\section{Proofs of the Mean Proposition and the Bridge}\label{sec:proof-bridge}
We first complete the mean proposition, then establish the higher-moment and excess-moment bridge inequalities.

\subsection{Exact Mean Formula and Asymptotics}\label{subsec:proof-mean}
Set $M=k+1$.
The mean identity $\E Z_{k,n} = \E Y_{M,n}$ was proved in Section~\ref{subsec:mean-identity-ideas}. We establish the exact formula and asymptotics from Proposition~\ref{prop:bridge} here. The results in this subsection largely build on existing work: the permutation representation for
the moments of $Y_{M,n}$ follows from \citet{dubach:2024}, while the mean asymptotics
for $Z_{k,n}$ were established by \citet{seleznova:2026}. We bring these results together to obtain the exact mean formula and the asymptotics needed for our subsequent analysis.

Following the notation from Section~\ref{sec:proof-ideas-higher-moments}, let $\tau_1 := (0_1,\dots, k_1)$ be a fixed cycle on $M$ elements of $I_1:=\{0_1,\dots k_1\}$ and let $\mathcal C_{M,1}$ denote the set
of all $M$-cycles on $I_1$.

\begin{proof}[Proof of Proposition~\ref{prop:bridge}]
Fix $M\geq1$, let $\cyc(\pi)$ denote the number of cycles of a permutation $\pi$, and define
\begin{align*}
    P_M(x):=\sum_{c\in\mathcal C_{M,1}} x^{\cyc(\tau_1^{-1}c)}.
\end{align*}
Lemma~9 of \citet{dubach:2024} connects $P_M$ with the moments of $Y_{M,n}$. In particular,
\begin{align*}
    n^M\E[Y_{M,n}]
    = P_M(n).
\end{align*}
Note that this permutation representation can also be obtained from Wick's formula, by an argument analogous to the proof of Lemma~\ref{lem:Z_commutator_representation}. 

For every $c\in\mathcal C_{M,1}$, the permutation $\tau_1^{-1}c$ acts on $M$ elements and hence has between one and $M$ cycles. Thus $P_M$ has degree at most $M$ and $P_M(0)=0$.

For $1\leq n\leq M$, applying the calculation in
\citet[eqs.~(2.6)--(2.8)]{dubach:2024} to the matrix $\sqrt n\,\W$ gives
\begin{align*}
M(M+1)P_M(n)=\prod_{j=0}^M(n+j).
\end{align*}
Thus $\prod_{j=0}^M(x+j)-M(M+1)P_M(x)$ is a monic polynomial
of degree $M+1$ with roots $0,1,\ldots,M$.
It therefore equals 
\begin{align*}
    \prod_{j=0}^M(x+j)-M(M+1)P_M(x) = \prod_{j=0}^M(x-j).
\end{align*}
Evaluating this identity at $x=n$ and using the mean identity $\E Z_{k,n}=\E Y_{M,n}$ yields the exact mean formula~\eqref{eq:mean_expression}, namely
\begin{equation}
    \E Z_{k,n}
    =
    \E Y_{M,n}
    =
    \frac{1}{M(M+1)n^M}
    \left[
        \prod_{j=0}^M(n+j)
        -
        \prod_{j=0}^M(n-j)
    \right].
\end{equation}
This identity holds for every $n\geq1$.

The asymptotic analysis of the mean expression~\eqref{eq:mean_expression} was carried out in Appendix~C of \citet{seleznova:2026}.
We briefly recall the main arguments here. For $M=o(n^{2/3})$, the starting point is the
Taylor expansion
\begin{align}\label{eq:taylor}
    \sum_{j=0}^M \log(1\pm j/n)
    = \pm\frac{M(M+1)}{2n}
      + O\!\left(\frac{M^3}{n^2}\right).
\end{align}
If $M/\sqrt n\to c\in(0,\infty)$, this yields
$\prod_{j=0}^M(1\pm j/n)\to e^{\pm c^2/2}$,
and hence
\[
    \E Y_{M,n}
    \longrightarrow
    \frac{2}{c^2}\sinh\!\left(\frac{c^2}{2}\right).
\]
Similarly, if $\sqrt n\ll M\ll n^{2/3}$,
the negative-sign product is negligible compared
with the positive-sign product, and
\eqref{eq:taylor} yields the scaling~\eqref{eq:mean-regimes2}:
\[
    \E Y_{M,n}
    \sim \frac{n}{M^2}
    \exp\left(\frac{M^2}{2n}\right)
    \sim \frac{n}{k^2}
    \exp\left(\frac{k^2}{2n}\right).
\]

The regime $M/\sqrt n\to0$, where cancellation
between the two products requires a sharper
estimate, is treated in
\citet[Theorem C.1]{seleznova:2026}.

Finally, the divergence when $M/\sqrt n\to\infty$
can be seen directly from the exact mean formula.
Indeed, expanding the difference of products,
all odd-degree terms are nonnegative.
Retaining only the cubic term yields
\[
    \E Y_{M,n}
    \geq
    \frac{2}{M(M+1)n^2}
    \sum_{1\leq i<j<\ell\leq M}ij\ell
    \geq
    C\frac{M^4}{n^2}
\]
for some absolute constant $C>0$ and all
sufficiently large $M$. This completes the proof of the scaling regimes in~\eqref{eq:mean-regimes} and  Proposition~\ref{prop:bridge}.
\end{proof}

\subsection{Permutation Representations for Higher Moments}\label{subsec:bridge-representations}

For a finite set $A$, write $S(A)$ for its permutation group,
and let $\cyc(\pi)$ denote the number of cycles of a
permutation $\pi$, including fixed points.
For $i=1,\ldots,p$, introduce the disjoint labeled sets
\begin{align*}
    I_i=\{0_i,\ldots,k_i\},
    \qquad
    I_i^\circ=\{1_i,\ldots,k_i\}.
\end{align*}
As in Section~\ref{sec:proof-ideas-higher-moments}, define the following fixed permutation
\begin{align*}
    \tau_p
    :=
    (0_1,1_1,\ldots,k_1)\cdots(0_p,1_p,\ldots,k_p)
    \in S\left(\bigcup_{i=1}^p I_i\right),
\end{align*}
let $\mathcal C_{M,p}$ denote the set of permutations of $\bigcup_{i=1}^p I_i$ consisting of $p$ cycles of length $M$, and further define 
\begin{align*}
    \mathcal C_{M,p}^{0}
    :=
    \left\{
        c\in\mathcal C_{M,p}:
        0_1,\ldots,0_p
        \text{ belong to distinct cycles of }c
    \right\}.
\end{align*}

We first give a permutation representation of $Y_{M,n}$ moments.
\begin{lemma}\label{lem:y-perm}
For all integers $n,M,p\geq1$, 
\begin{align}
    \mathbb E[Y_{M,n}^p]
    &=
    \frac{1}{M^p n^{pM}}
    \sum_{\sigma\in S(\cup_{i=1}^p I_i)}
    n^{\cyc(\tau_p^{-1}\sigma\tau_p\sigma^{-1})}
    \notag\\
    &=
    \frac{p!}{n^{pM}}
    \sum_{c\in\mathcal C_{M,p}}
    n^{\cyc(\tau_p^{-1}c)}.
    \label{eq:y-perm}
\end{align}
\end{lemma}

\begin{proof}
The first identity follows from Lemma~9 of
\citet{dubach:2024}. One can also obtain it via Wick expansion as in the proof of the next Lemma~\ref{lem:Z_commutator_representation}.

For the second identity, consider the map
\[
    \Psi_{\tau_p}:
    S\left(\bigcup_{i=1}^p I_i\right)
    \longrightarrow \mathcal C_{M,p},
    \qquad
    \sigma\longmapsto\sigma\tau_p\sigma^{-1}.
\]
Fix $c\in\mathcal C_{M,p}$ and label its $p$ cycles as
$C_1,\ldots,C_p$. The condition $\Psi_{\tau_p}(\sigma)=c$
is equivalent to
\[
    \sigma\tau_p=c\sigma.
\]
Thus, $\sigma$ maps each cycle $I_i$ of $\tau_p$ onto a cycle of $c$, preserving the cyclic order. Since $\sigma$ is a bijection,
distinct cycles $I_i$ must map onto distinct cycles of $c$. There are therefore $p!$ choices for the assignment of the source cycles $I_1,\ldots,I_p$ to the target cycles $C_1,\ldots,C_p$.

For a fixed assignment, choose the image $u_i:=\sigma(0_i)$ in the target cycle assigned to $I_i$. There are $M$ choices
for each $u_i$. Once $u_i$ is chosen, the relation $\sigma\tau_p=c\sigma$ determines all remaining images on $I_i$:
\[
    \sigma(r_i)=c^r(u_i),\qquad r=0,\ldots,M-1.
\]
Conversely, every such assignment and choice of starting points defines a permutation $\sigma$ satisfying $\sigma\tau_p=c\sigma$.
Different choices give different permutations, and every preimage
arises in this way. Hence $\Psi_{\tau_p}$ is surjective and
\[
    |\Psi_{\tau_p}^{-1}(c)|=p!M^p
    \qquad\text{for every }c\in\mathcal C_{M,p}.
\]
Finally, the summand depends on $\sigma$ only through
$c=\Psi_{\tau_p}(\sigma)$. Grouping the terms by these preimages gives
\[
    \sum_{\sigma\in S(\cup_{i=1}^p I_i)}
        n^{\cyc(\tau_p^{-1}\sigma\tau_p\sigma^{-1})}
    =p!M^p\sum_{c\in\mathcal C_{M,p}}
        n^{\cyc(\tau_p^{-1}c)}.
\]
Dividing by $M^p n^{pM}$ yields the second identity.
\end{proof}

We next give a corresponding representation for the moments of $Z_{k,n}$ as a sum over permutations.

\begin{lemma}\label{lem:Z_commutator_representation}
For $\sigma\in S(\bigcup_{i=1}^p I_i^\circ)$, denote by
$\widetilde\sigma\in S(\bigcup_{i=1}^p I_i)$ its extension defined by
\[
    \widetilde\sigma(0_i)=0_i,
    \qquad i=1,\ldots,p.
\] Then, for all integers
$n,k,p\geq1$, with $M=k+1$,
\begin{align}
    \E[Z_{k,n}^p]
    &=
    \frac{1}{n^{pM}}
    \sum_{\sigma\in S(\cup_{i=1}^p I_i^\circ)}
    n^{\cyc(
        \tau_p^{-1}\widetilde\sigma^{-1}
        \tau_p\widetilde\sigma
    )}
    \notag\\
    &=
    \frac{1}{n^{pM}}
    \sum_{c\in\mathcal C_{M,p}^{0}}
    n^{\cyc(\tau_p^{-1}c)}.
    \label{eq:z-perm}
\end{align}
The second equality follows from the fact that the map
\[
    \Phi_{\tau_p}:
    S\left(\bigcup_{i=1}^p I_i^\circ\right)
    \longrightarrow
    \mathcal C_{M,p}^{0},
    \qquad
    \sigma\longmapsto
    \widetilde\sigma^{-1}\tau_p\widetilde\sigma,
\]
is a bijection.
\end{lemma}

\begin{proof}
\noindent\textbf{Part 1: Wick expansion.}
We start by expanding the $p$ copies of the trace:
\begin{align}
    \tr^p\bigl((\W^k)^*\W^k\bigr)
    &=
    \sum_{\substack{
        a^{(i)}_0,\dots,a^{(i)}_k=1\\
        b^{(i)}_1,\dots,b^{(i)}_{k-1}=1\\
        i=1,\dots,p
    }}^n
    \prod_{i=1}^p
    \overline{\W}_{a^{(i)}_0a^{(i)}_1}
    \cdots
    \overline{\W}_{a^{(i)}_{k-1}a^{(i)}_k}
    \W_{a^{(i)}_0b^{(i)}_1}
    \cdots
    \W_{b^{(i)}_{k-1}a^{(i)}_k}.
    \label{eq:trace_p_expansion}
\end{align}
We will rewrite the expectation of the above monomials in complex Gaussian entries of $\W$ using the Wick/Isserlis theorem~\citep{isserlis1918formula}. For proper complex Gaussians, any Wick pairing that contains
a pair of two conjugate or two non-conjugate entries gives zero
contribution. Hence, the non-zero Wick pairings are indexed by
permutations $\sigma\in S(I^\circ)$, where $I_i^\circ=\{1_i,\dots,k_i\}$ and
$I^\circ:=\cup_{i=1}^p I_i^\circ$. For $r_i\in I^\circ$, write $\sigma(r_i)=s_j$. 
Moreover, we set
\begin{equation}
     a^{(i)}_{0}=b^{(i)}_{0},
    \qquad
    a^{(i)}_k=b^{(i)}_k, \quad i=1,\dots, p.
    \label{eq:wick_boundary_conditions}
\end{equation}
We can then write the
Wick expansion as
\begin{align}
    \E\left[
        \tr^p\bigl((\W^k)^*\W^k\bigr)
    \right]
    &=
    \sum_{\sigma\in S(I^\circ)}
    \sum_{\substack{
        a^{(i)}_0,\dots,a^{(i)}_k=1\\
        b^{(i)}_1,\dots,b^{(i)}_{k-1}=1\\
        i=1,\dots,p
    }}^n
    \prod_{\substack{r_i\in I^\circ\\
                     \sigma(r_i)=s_j}}
    \E\left[
        \overline{\W}_{a^{(i)}_{r-1}a^{(i)}_r}
        \W_{b^{(j)}_{s-1}b^{(j)}_s}
    \right].
    \label{eq:wick_p_expansion}
\end{align}
Since the entries of $\W$ have variance $1/n$, we have
\begin{equation}
    \E\left[
        \overline{\W}_{a^{(i)}_{r-1}a^{(i)}_r}
        \W_{b^{(j)}_{s-1}b^{(j)}_s}
    \right]
    =
    \frac{1}{n}
    \delta_{a^{(i)}_{r-1},b^{(j)}_{s-1}}
    \delta_{a^{(i)}_r,b^{(j)}_s}.
    \label{eq:wick_delta_p}
\end{equation}
Thus, a permutation $\sigma$ with $\sigma(r_i)=s_j$ induces the index
equalities
\begin{equation}
    a^{(i)}_{r-1}=b^{(j)}_{s-1},
    \qquad
    a^{(i)}_r=b^{(j)}_s, \quad r_i\in I^\circ.
    \label{eq:wick_identifications_p}
\end{equation}
Let $F_p(\sigma)$ denote the number of free index equivalence classes
left by these identifications. Since there are $pk$ Wick pairs, each
contributing a factor $1/n$, we obtain
\begin{equation}
    \E\left[
        \tr^p\bigl((\W^k)^*\W^k\bigr)
    \right]
    =
    \frac{1}{n^{pk}}
    \sum_{\sigma\in S(I^\circ)}
    n^{F_p(\sigma)}.
    \label{eq:moment_free_indices_p}
\end{equation}
It remains to identify $F_p(\sigma)$.

\noindent\textbf{Part 2: Commutator representation of $F_p(\sigma)$.} Eliminating the $b$-indices between two consecutive relations in~\eqref{eq:wick_identifications_p} gives,
for $s=1,\dots,k-1$,
\begin{equation}
    a^{(i)}_r=a^{(i')}_{r'-1},
    \qquad
    \text{where }
    \sigma^{-1}(s_j)=r_i,
    \quad
    \sigma^{-1}((s+1)_j)=r'_{i'}.
    \label{eq:a_cond_p}
\end{equation}
At the two boundaries, the convention $b^{(j)}_0=a^{(j)}_0$ and
$b^{(j)}_k=a^{(j)}_k$ gives
\begin{equation}
    a^{(i)}_{r-1}=a^{(j)}_0,
    \qquad
    \text{where }\sigma^{-1}(1_j)=r_i,
    \label{eq:a_cond_boundary_left_p}
\end{equation}
and
\begin{equation}
    a^{(i)}_r=a^{(j)}_k,
    \qquad
    \text{where }\sigma^{-1}(k_j)=r_i.
    \label{eq:a_cond_boundary_right_p}
\end{equation}

We now use the zero points to write these identifications in a more
convenient form. Recall that $I:=\cup_{i=1}^p I_i$ and
\begin{equation}
    \tau_p
    =
    (0_1,1_1,\dots,k_1)\cdots
    (0_p,1_p,\dots,k_p),
\end{equation}
while $\widetilde\sigma$ denotes the extension of $\sigma$ to $I$
which fixes $0_1,\dots,0_p$.

For $r_i\in I_i$, write $\mathfrak a_{r_i}:=a^{(i)}_r$. Combining
\eqref{eq:a_cond_p}--\eqref{eq:a_cond_boundary_right_p}, we obtain
\begin{equation}
    \mathfrak a_u
    =
    \mathfrak a_{
        \tau_p^{-1}\widetilde\sigma^{-1}
        \tau_p\widetilde\sigma(u)
    },
    \qquad
    u\in I.
    \label{eq:a_identification_p}
\end{equation}
Indeed, the relations \eqref{eq:a_cond_p} give the interior cases,
while \eqref{eq:a_cond_boundary_left_p} and
\eqref{eq:a_cond_boundary_right_p} give the corresponding
identifications involving the zero points.

Let
\begin{equation}
    g_\sigma
    :=
    \tau_p^{-1}\widetilde\sigma^{-1}
    \tau_p\widetilde\sigma.
    \label{eq:def_g_sigma_p}
\end{equation}
Then \eqref{eq:a_identification_p} simply says that
$\mathfrak a_u=\mathfrak a_{g_\sigma(u)}$ for every $u\in I$.

Now consider a cycle $(u_1,\dots,u_m)$ of $g_\sigma$. The previous
identity implies that
$\mathfrak a_{u_1}=\cdots=\mathfrak a_{u_m}$, so all positions in the
same cycle of $g_\sigma$ belong to the same equivalence class.
Conversely, \eqref{eq:a_identification_p} only identifies a position
$u$ with $g_\sigma(u)$ and therefore cannot identify positions
belonging to different cycles. Since all $b$-indices have already
been eliminated, there are no further index identifications.
Consequently, the equivalence classes of free indices are precisely
the cycles of $g_\sigma$, and therefore
\begin{equation}
    F_p(\sigma)
    =
    \cyc\left(
        \tau_p^{-1}\widetilde\sigma^{-1}
        \tau_p\widetilde\sigma
    \right).
    \label{eq:Fp_cycles}
\end{equation}

Substituting \eqref{eq:Fp_cycles} into
\eqref{eq:moment_free_indices_p} gives
\begin{equation}
    \E\left[
        \tr^p\bigl((\W^k)^*\W^k\bigr)
    \right]
    =
    \frac{1}{n^{pk}}
    \sum_{\sigma\in S(\cup_{i=1}^p I_i^\circ)}
    n^{\cyc(
        \tau_p^{-1}\widetilde\sigma^{-1}
        \tau_p\widetilde\sigma
    )}.
\end{equation}
Finally, using $Z_{k,n}=n^{-1}\tr((\W^k)^*\W^k)$ proves
\begin{equation}
    \E[Z_{k,n}^p]
    =
    \frac{1}{n^{p(k+1)}}
    \sum_{\sigma\in S(\cup_{i=1}^p I_i^\circ)}
    n^{\cyc(
        \tau_p^{-1}\widetilde\sigma^{-1}
        \tau_p\widetilde\sigma
    )}.
\end{equation}
\medskip
\noindent\textbf{Part 3: The bijection $\Phi_{\tau_p}$.}
First, we check that $\Phi_{\tau_p}$ indeed takes values in
$\mathcal C_{k+1,p}^{0}$. For $\sigma\in S(\cup_{i=1}^p I_i^\circ)$,
the permutation
\begin{equation}
    c
    :=
    \Phi_{\tau_p}(\sigma)
    =
    \widetilde\sigma^{-1}\tau_p\widetilde\sigma
\end{equation}
is conjugate to $\tau_p$. Since $\tau_p$ consists of $p$ cycles of
length $k+1$, the same is true of $c$. Moreover,
$\widetilde\sigma(0_i)=0_i$ for every $i=1,\dots,p$. Therefore, the
cycle of $c$ containing $0_i$ is
\begin{equation}
    \bigl(
        \widetilde\sigma^{-1}(0_i),\,
        \widetilde\sigma^{-1}(1_i),\,
        \dots,\,
        \widetilde\sigma^{-1}(k_i)
    \bigr) = \bigl(
        0_i,\,
        \sigma^{-1}(1_i),\,
        \dots,\,
        \sigma^{-1}(k_i)
    \bigr).
    \label{eq:cycles_conjugated_tau}
\end{equation}
In particular, $0_1,\dots,0_p$ belong to distinct cycles of $c$, and
hence $c\in\mathcal C_{k+1,p}^{0}$.

To prove injectivity, suppose that for some
$\sigma_1,\sigma_2\in S(\cup_{i=1}^p I_i^\circ)$,
\begin{equation}
\widetilde\sigma_1^{-1}\tau_p\widetilde\sigma_1
    =
\widetilde\sigma_2^{-1}\tau_p\widetilde\sigma_2.
\end{equation}
Then
$\widetilde\sigma_2\widetilde\sigma_1^{-1}\tau_p
=
\tau_p\widetilde\sigma_2\widetilde\sigma_1^{-1}$, so
$g:=\widetilde\sigma_2\widetilde\sigma_1^{-1}$ belongs to the
centralizer of $\tau_p$.
The centralizer of $\tau_p$ consists of independent rotations of its
$p$ cycles, together with permutations of cycles of the same length.
On the other hand, $g$ fixes each of the points $0_1,\dots,0_p$
individually. Hence, $g$ cannot permute the cycles of $\tau_p$, and
none of the cycles can be rotated nontrivially. Therefore $g=Id$,
and thus $\widetilde\sigma_1=\widetilde\sigma_2$, which implies
$\sigma_1=\sigma_2$. Hence, $\Phi_{\tau_p}$ is injective.

It remains to show that the domain and $\mathcal C_{k+1,p}^{0}$ have
the same cardinality. Clearly,
\begin{equation}
    \left|S\left(\bigcup_{i=1}^p I_i^\circ\right)\right|
    =
    (pk)!.
\end{equation}
To construct an element of $\mathcal C_{k+1,p}^{0}$, we first
partition the $pk$ points in $\cup_{i=1}^p I_i^\circ$ into $p$
sets of size $k$, where the $i$-th set contains the points that accompany $0_i$ in its cycle. We then choose the order of the $k$
points within each cycle. Therefore,
\begin{equation}
    |\mathcal C_{k+1,p}^{0}|
    =
    \binom{pk}{k,\dots,k}(k!)^p
    =
    (pk)!.
\end{equation}
Thus, $\Phi_{\tau_p}$ is an injective map between two finite sets of
the same cardinality, and is therefore a bijection.
\end{proof}
\subsection{Comparison of Weighted Permutation Sums}\label{subsec:bridge-comparison}
We now compare the weighted sums over $\mathcal C_{M,p}^{0}$
and $\mathcal C_{M,p}$ to obtain the desired moment bounds.
\begin{lemma}\label{lem:cycle_sum_comparison}
Let $\tau_p:=(0_1,1_1,\ldots,k_1)\cdots(0_p,1_p,\ldots,k_p)$,
and let $\mathcal C_{k+1,p}$ and $\mathcal C_{k+1,p}^{0}$ be as defined above. Then, for every $n,k, p$,
\begin{equation}\label{eq:cycle_sum_comparison}
    \frac{p!}{p^p}
    \sum_{c\in\mathcal C_{k+1,p}}
        n^{\cyc(\tau_p^{-1}c)}
    \leq
    \sum_{c\in\mathcal C_{k+1,p}^{0}}
        n^{\cyc(\tau_p^{-1}c)}
    \leq
    \sum_{c\in\mathcal C_{k+1,p}}
        n^{\cyc(\tau_p^{-1}c)}.
\end{equation}
\end{lemma}

\begin{proof}[Proof of Lemma~\ref{lem:cycle_sum_comparison}]
The upper bound follows immediately from
$\mathcal C_{k+1,p}^{0}\subseteq\mathcal C_{k+1,p}$, so it remains
to prove the lower bound.

For $\bu=(u_1,\dots,u_p)\in I_1\times\cdots\times I_p$, let
$\mathcal C_{k+1,p}^{\bu}$ denote the set of permutations
$c\in\mathcal C_{k+1,p}$ for which the points $u_1,\dots,u_p$
belong to distinct cycles of $c$. Thus,
$\mathcal C_{k+1,p}^{(0_1,\dots,0_p)}
=\mathcal C_{k+1,p}^{0}$.

We first notice that the sum does not depend
on the choice of $\bu$. Indeed, for each $i=1,\dots,p$, choose
$m_i\in\{0,\dots,k\}$ such that $u_i=\tau_p^{m_i}(0_i)$, and define
$\rho_{\bu}$ by letting it act as $\tau_p^{m_i}$ on $I_i$. Then
$\rho_{\bu}$ commutes with $\tau_p$ and maps $0_i$ to $u_i$ for
every $i$. Conjugation by $\rho_{\bu}$ therefore gives a bijection
from $\mathcal C_{k+1,p}^{0}$ to
$\mathcal C_{k+1,p}^{\bu}$. Moreover,
\begin{equation}
    \cyc\left(
        \tau_p^{-1}\rho_{\bu}c\rho_{\bu}^{-1}
    \right)
    =
    \cyc\left(
        \rho_{\bu}\tau_p^{-1}c\rho_{\bu}^{-1}
    \right)
    =
    \cyc(\tau_p^{-1}c).
\end{equation}
Consequently,
\begin{equation}\label{eq:marked_cycle_sum}
    \sum_{c\in\mathcal C_{k+1,p}^{\bu}}
        n^{\cyc(\tau_p^{-1}c)}
    =
    \sum_{c\in\mathcal C_{k+1,p}^{0}}
        n^{\cyc(\tau_p^{-1}c)}
\end{equation}
for every $\bu\in I_1\times\cdots\times I_p$.

We now average \eqref{eq:marked_cycle_sum} over all
$(k+1)^p$ choices of $\bu$. This gives
\begin{equation}\label{eq:average_marked_cycle_sum}
    \sum_{c\in\mathcal C_{k+1,p}^{0}}
        n^{\cyc(\tau_p^{-1}c)}
    =
    \sum_{c\in\mathcal C_{k+1,p}}
        w(c)\,n^{\cyc(\tau_p^{-1}c)},
\end{equation}
where $w(c)$ is the proportion of choices
$\bu\in I_1\times\cdots\times I_p$ for which
$u_1,\dots,u_p$ belong to distinct cycles of $c$.

It remains to bound $w(c)$. Fix $c\in\mathcal C_{k+1,p}$ and denote
its cycles by $C_1,\dots,C_p$. Define the $p\times p$ matrix
$D(c)=(d_{ij})$ by
\begin{equation}
    d_{ij}:=|I_i\cap C_j|.
\end{equation}
Since each $I_i$ and each $C_j$ contains exactly $k+1$ points, every
row and every column of $D(c)$ sums to $k+1$. Hence,
$D(c)/(k+1)$ is doubly stochastic.

The number of choices $\bu=(u_1,\dots,u_p)$ for which
$u_1,\dots,u_p$ lie in distinct cycles of $c$ is
\begin{equation}
    \sum_{\pi\in S_p}
        \prod_{i=1}^p d_{i,\pi(i)}
    =
    \operatorname{perm}D(c).
\end{equation}
Therefore,
\begin{equation}\label{eq:w_permanent}
    w(c)
    =
    \frac{\operatorname{perm}D(c)}{(k+1)^p}
    =
    \operatorname{perm}\left(\frac{D(c)}{k+1}\right).
\end{equation}
By the van der Waerden inequality for permanents of doubly stochastic
matrices~\citep{egorychev1981proof},
\begin{equation}
    w(c)\geq \frac{p!}{p^p}.
\end{equation}
On the other hand, $w(c)\leq 1$ since it is a proportion.

Substituting these bounds into
\eqref{eq:average_marked_cycle_sum} gives
\begin{equation}
    \frac{p!}{p^p}
    \sum_{c\in\mathcal C_{k+1,p}}
        n^{\cyc(\tau_p^{-1}c)}
    \leq
    \sum_{c\in\mathcal C_{k+1,p}^{0}}
        n^{\cyc(\tau_p^{-1}c)}
    \leq
    \sum_{c\in\mathcal C_{k+1,p}}
        n^{\cyc(\tau_p^{-1}c)},
\end{equation}
which proves the claim.
\end{proof}

Combining Lemmas~\ref{lem:Z_commutator_representation} %
and~\ref{lem:cycle_sum_comparison}, we obtain the first assertion \eqref{eq:thm:bridge-variance1} of Theorem \ref{thm:bridge-variance}.
For $p=1$, this recovers the mean identity already proven in Section~\ref{subsec:mean-identity-ideas}.
These results also give
\[
    R^{(p)}_{k,n}
    \leq \frac{1}{p!}q^{(p)}_{k+1,n},
\]
which is the upper bound in the second assertion \eqref{eq:thm:bridge-variance2}.
\subsection{Bridge for Excess Moments}\label{subsec:bridge-excess}
To establish the lower bound in \eqref{eq:thm:bridge-variance2}, that is,
$
    \frac{1}{p^p}q^{(p)}_{k+1,n}\leq R^{(p)}_{k,n},
$
we use the following corollary of
Lemma~\ref{lem:cycle_sum_comparison}.
\begin{corollary}\label{cor:connected_cycle_sum_comparison}
Let $\mathcal C_{k+1,p}^{\mathrm{disc}}$ denote the set of $c\in\mathcal C_{k+1,p}$ whose cycle supports are exactly $I_1,\dots,I_p$, and let
\begin{equation}
    \mathcal C_{k+1,p}^{\mathrm{conn}}
    :=
    \mathcal C_{k+1,p}
    \setminus
    \mathcal C_{k+1,p}^{\mathrm{disc}}.
\end{equation}
Here ``connected'' means the complement of the completely blockwise-disconnected contribution; partially disconnected configurations are included. We define $\mathcal C_{k+1,p}^{0,\mathrm{disc}}$ and
$\mathcal C_{k+1,p}^{0,\mathrm{conn}}$ similarly. Then
\begin{equation}\label{eq:connected_cycle_sum_comparison}
    \frac{p!}{p^p}
    \sum_{c\in\mathcal C_{k+1,p}^{\mathrm{conn}}}
        n^{\cyc(\tau_p^{-1}c)}
    \leq
    \sum_{c\in\mathcal C_{k+1,p}^{0,\mathrm{conn}}}
        n^{\cyc(\tau_p^{-1}c)}
    \leq
    \sum_{c\in\mathcal C_{k+1,p}^{\mathrm{conn}}}
        n^{\cyc(\tau_p^{-1}c)}.
\end{equation}
\end{corollary}

\begin{proof}[Proof of Corollary~\ref{cor:connected_cycle_sum_comparison}]
Notice first that every
$c\in\mathcal C_{k+1,p}^{\mathrm{disc}}$ also belongs to
$\mathcal C_{k+1,p}^{0}$, since its cycle supports are exactly
$I_1,\dots,I_p$. Therefore,
\begin{equation}
    \mathcal C_{k+1,p}^{0,\mathrm{disc}}
    =
    \mathcal C_{k+1,p}^{\mathrm{disc}}.
\end{equation}
Moreover, for such a $c$, any choice of one point from each $I_i$
puts the chosen points in distinct cycles of $c$. Hence the weights defined in~\eqref{eq:average_marked_cycle_sum} satisfy $w(c)=1$ for
every $c\in\mathcal C_{k+1,p}^{\mathrm{disc}}$ and we can subtract the disconnected part from both sides of~\eqref{eq:average_marked_cycle_sum} to get
\begin{equation}\label{eq:connected_weighted_identity}
    \sum_{c\in\mathcal C_{k+1,p}^{0,\mathrm{conn}}}
        n^{\cyc(\tau_p^{-1}c)}
    =\sum_{c\in\mathcal C_{k+1,p}^{\mathrm{conn}}}
        w(c)n^{\cyc(\tau_p^{-1}c)}.
\end{equation}
The result now follows directly from
$p!/p^p\leq w(c)\leq1$,  which is proven in Lemma~\ref{lem:cycle_sum_comparison}.
\end{proof}

For $c\in\mathcal C_{M,p}^{\mathrm{disc}}$, the permutation
$\tau_p^{-1}c$ acts separately on the blocks $I_1,\ldots,I_p$.
The corresponding sum 
$$
\sum_{c\in\mathcal C_{M,p}^{\mathrm{disc}}}
n^{\cyc(\tau_p^{-1}c)}
$$
therefore factors into a product of $p$ identical sums,
each equal to $n^M\mathbb E[Z_{k,n}]$ by
Lemma~\ref{lem:Z_commutator_representation}
with $p=1$. That is, \[
\sum_{c\in\mathcal C_{M,p}^{\mathrm{disc}}}
n^{\cyc(\tau_p^{-1}c)}
=
\bigl(n^M\mathbb E[Z_{k,n}]\bigr)^p
=
n^{pM}\bigl(\mathbb E[Z_{k,n}]\bigr)^p.
\]
Moreover,
$\mathcal C_{M,p}^{0,\mathrm{disc}}
=\mathcal C_{M,p}^{\mathrm{disc}}$.
Subtracting the $(\mathbb E[Z_{k,n}])^p$ 
 from the representation of $\mathbb E[Z_{k,n}^p]$
in Lemma~\ref{lem:Z_commutator_representation}, we get
\begin{align} \label{eq:n1}
    \mathbb E[Z_{k,n}^p]
    -\bigl(\mathbb E[Z_{k,n}]\bigr)^p
    &=
    \frac{1}{n^{pM}}
    \sum_{c\in\mathcal C_{M,p}^{0,\mathrm{conn}}}
        n^{\cyc(\tau_p^{-1}c)}.
\end{align}
Similarly, using \eqref{eq:y-perm} and  $\mathbb E[Z_{k,n}]=\mathbb E[Y_{M,n}]$,
we also get 
\begin{align} \label{eq:n2}
     \mathbb E[Y_{M,n}^p]
    -p!\bigl(\mathbb E[Y_{M,n}]\bigr)^p
    &=
    \frac{p!}{n^{pM}}
    \sum_{c\in\mathcal C_{M,p}^{\mathrm{conn}}}
        n^{\cyc(\tau_p^{-1}c)}.
\end{align}
Applying Corollary~\ref{cor:connected_cycle_sum_comparison}
to \eqref{eq:n1} and \eqref{eq:n2} yields
\begin{align*}
    \frac{1}{p^p}
    \left(
        \mathbb E[Y_{M,n}^p]
        -p!\bigl(\mathbb E[Y_{M,n}]\bigr)^p
    \right)
    &\leq
    \mathbb E[Z_{k,n}^p]
    -\bigl(\mathbb E[Z_{k,n}]\bigr)^p
    \\ &\leq
    \frac{1}{p!}
    \left(
        \mathbb E[Y_{M,n}^p]
        -p!\bigl(\mathbb E[Y_{M,n}]\bigr)^p
    \right).
\end{align*}
Dividing by $(\mathbb E[Y_{M,n}])^p$ and using
Proposition~\ref{prop:bridge}, we obtain
\begin{align*}
    \frac{1}{p^p}
    \left(
        \frac{\mathbb E[Y_{M,n}^p]}
             {(\mathbb E[Y_{M,n}])^p}
        -p!
    \right)
    &\leq
    \frac{\mathbb E[Z_{k,n}^p]}
         {(\mathbb E[Z_{k,n}])^p}
    -1
    \leq
    \frac{1}{p!}
    \left(
        \frac{\mathbb E[Y_{M,n}^p]}
             {(\mathbb E[Y_{M,n}])^p}
        -p!
    \right).
\end{align*}
This proves the remaining lower bound in
\eqref{eq:thm:bridge-variance2} and completes the proof of
Theorem~\ref{thm:bridge-variance}.

\section{Proofs of the Scaling Regimes}
\label{sec:proof-q-asymptotics}

Throughout this section, fix an integer $p\geq2$ and set $M=k+1$.
For convenience, we recall the definition of $s_{p,n}$ from~\eqref{eq:fluc_scaling}:
\[
    s_{p,n}
    =\frac{k^2}{n}-\frac{1}{p}\log n+\frac{1}{p}\log\log n.
\]
We first reduce excess moments to normalized cumulants and state the radial dominance theorem. We then prove the intermediate-power limits, followed by the small- and large-power regimes. The three conclusions are given below. Together, they prove Theorem~\ref{thm:q-asymptotics}.

In Section~\ref{subsec:q-small-powers}, we consider
$M=O(\sqrt n)$ and prove that
\begin{align}
    q^{(p)}_{M,n}\longrightarrow0.
    \label{eq:q-small-powers}
\end{align}
In Section~\ref{subsec:q-intermediate-powers}, we consider
$\sqrt n\ll M\ll n^{2/3}$ and prove that
\begin{align}
    q^{(p)}_{M,n}
    \longrightarrow
    \begin{cases}
        0,
        & s_{p,n}\to-\infty,\\[1mm]
        p^{-(p+1)/2}
        \exp\!\left(\dfrac{p(p-1)}{2}s\right),
        & s_{p,n}\to s\in\mathbb R,\\[1mm]
        \infty,
        & s_{p,n}\to+\infty.
    \end{cases}
    \label{eq:q-intermediate-powers}
\end{align}
Finally, in Section~\ref{subsec:q-large-powers}, we consider
$M\geq c n^{2/3}$ for any fixed $c>0$ and prove that
\begin{align}
    q^{(p)}_{M,n}\longrightarrow\infty.
    \label{eq:q-large-powers}
\end{align}

\subsection{Reduction to Normalized Cumulants}\label{subsec:q-auxiliary}
We now reduce the excess moments to normalized mixed cumulants. Set
\begin{align} \label{eq:def_K}
    T:=\tr(\mathbf W^M),
    \qquad
    K_q:=\kappa(
        \underbrace{T,\ldots,T}_{q},
        \underbrace{\overline T,\ldots,\overline T}_{q}),
    \qquad q\geq1,
\end{align}
where $\kappa$ denotes the joint cumulant. 
In the following, we reduce the analysis of $q^{(p)}_{M,n}$
to that of the normalized cumulants $K_q/K_1^q$.
Since $Y_{M,n}=|T|^2/M$, we have
$K_1=\mathbb E|T|^2=M\mathbb E[Y_{M,n}]$.
By rotational invariance, mixed cumulants containing
unequal numbers of $T$ and $\overline T$ vanish.
The moment-cumulant formula therefore gives
\begin{equation}\label{eq:moment-cumulant}
    p!+q^{(p)}_{M,n}
    =
    (p!)^2
    \sum_{\substack{m_1,\ldots,m_p\geq0\\
                    \sum_{q=1}^p qm_q=p}}
    \prod_{q=1}^p
    \frac{(K_q/K_1^q)^{m_q}}
         {(q!)^{2m_q}m_q!}.
\end{equation}
Here $m_q$ counts the number of blocks containing $q$
copies of $T$ and $q$ copies of $\overline T$.
The term with $m_1=p$ equals $p!$, while the term with
$m_p=1$ equals $K_p/K_1^p$.
Every term other than $m_1=p$ contains at least one
factor $K_q/K_1^q$ with $q\geq2$.

\subsection{Radial Dominance}\label{subsec:radial-statement}
The following result controls the cumulants through the radial moments of the eigenvalues of $\mathbf W$. Let $\lambda_1,\ldots,\lambda_n$ denote the eigenvalues of $\mathbf W$, and set
\begin{equation} \label{eq:def_A}
    A_L(n):=\mathbb E\left[\sum_{j=1}^n|\lambda_j|^{2L}\right],
    \qquad L\geq0.
\end{equation}
\begin{theorem}[Radial dominance]
\label{thm:radial-dominance}

Fix an integer $q\geq2$.
\begin{enumerate}
    \item For all integers $n\geq1$ and $L\geq0$,
    \begin{equation}
        A_L(n)=\frac{(n+L)!}{(L+1)(n-1)!\,n^L}.
        \label{eq:radial-term-exact}
    \end{equation}

    \item There exists a constant $C_q<\infty$ such that,
    for all integers $n,M\geq1$,
    \begin{equation}
        |K_q-A_{qM}(n)|
        \leq C_q A_{qM}(n)
        \exp\left(-\frac{(q-1)M^2}{n+qM}\right).
        \label{eq:radial-error-bound}
    \end{equation}
    In particular, $|K_q|\leq(1+C_q)A_{qM}(n)$.

    \item If $n\to\infty$ and $M/\sqrt n\to\infty$, then
    \[
        K_q\sim A_{qM}(n).
    \]
    If, in addition, $M=o(n^{2/3})$, then
    \begin{equation}
        K_q\sim\frac{n}{qM}
        \exp\left(\frac{q^2M^2}{2n}\right).
        \label{eq:radial-dominance-asymptotic}
    \end{equation}
\end{enumerate}
\end{theorem}
The proof of Theorem \ref{thm:radial-dominance} is deferred to Section \ref{sec:proof_thm_radial}. 

\subsection{\texorpdfstring{Intermediate Powers: Proof of \eqref{eq:q-intermediate-powers}}{Intermediate Powers}}
\label{subsec:q-intermediate-powers}
\begin{proof}[Proof of \eqref{eq:q-intermediate-powers}]
By \eqref{eq:radial-dominance-asymptotic} of Theorem~\ref{thm:radial-dominance} and
Proposition~\ref{prop:bridge}, for every $2\leq q\leq p$,
\begin{align}
    \frac{K_q}{K_1^q}
    &\sim
    \frac{1}{q}\left(\frac{M}{n}\right)^{q-1}
    \exp\left(\frac{q(q-1)M^2}{2n}\right).
    \label{eq:normalized-cumulants-intermediate}
\end{align}
It remains to determine the regime in which the right-hand side is of
order one. Set $x=M^2/n$, so that $M/n=(x/n)^{1/2}$. Then
\begin{equation}
    \frac{K_p}{K_1^p}
    \sim
    \frac{1}{p}
    \left(\frac{x}{n}\right)^{(p-1)/2}
    e^{p(p-1)x/2}.
\end{equation}
Thus $K_p/K_1^p$ is of order one when
$px+\log x-\log n=O(1)$. This condition implies $x\sim p^{-1}\log n$, which gives
$\log x=\log\log n+O(1)$, and therefore the scaling regime of interest is
\begin{equation}
    \frac{M^2}{n}
    =
    \frac{\log n-\log\log n}{p}
    +O(1).
\end{equation}
Since $M=k+1$, by definition of $s_{p,n}$~\eqref{eq:fluc_scaling}\footnote{Note that the $s_{p,n}$ is defined in terms of $k$ rather than $M$, but the change from $M^2/n$ to $k^2/n$ does not affect the resulting asymptotics. }, we may rewrite
\eqref{eq:normalized-cumulants-intermediate} as
\begin{equation}\label{eq:cumulants-spn}
    \frac{K_q}{K_1^q}
    \sim
    \frac{1}{q}
    \left[
        \frac{x}{\log n}
        \exp\left(p s_{p,n}-(p-q)x\right)
    \right]^{(q-1)/2},
    \qquad 2\leq q\leq p.
\end{equation}
Recall that $x=M^2/n\to\infty$.

Suppose first that $s_{p,n}\to s\in\mathbb R$.
Then $x/\log n\to1/p$, and
\eqref{eq:cumulants-spn} with $q=p$ gives
\begin{equation}
    \frac{K_p}{K_1^p}
    \longrightarrow
    p^{-(p+1)/2}
    \exp\left(\frac{p(p-1)}{2}s\right).
\end{equation}
For every $2\leq q<p$, \eqref{eq:cumulants-spn} also yields
$K_q/K_1^q\to0$, since $(p-q)x\to\infty$.
In the moment-cumulant formula
\eqref{eq:moment-cumulant}, every term other than
$m_1=p$ and $m_p=1$ contains at least one such
lower-order normalized cumulant and therefore vanishes.
Consequently,
\begin{equation}
    q^{(p)}_{M,n}
    \sim \frac{K_p}{K_1^p}
    \longrightarrow
    p^{-(p+1)/2}
    \exp\left(\frac{p(p-1)}{2}s\right).
\end{equation}

Next, suppose that $s_{p,n}\to-\infty$.
Since $s_{p,n}\leq0$ for all sufficiently large $n$,
its definition in \eqref{eq:fluc_scaling} and $\delta_n=o(1)$ imply $x/\log n=O(1)$.

Since $p-q\geq0$, \eqref{eq:cumulants-spn} implies
\[
    \frac{K_q}{K_1^q}\longrightarrow0,
    \qquad 2\leq q\leq p.
\]
Every term in \eqref{eq:moment-cumulant} other than
$m_1=p$ therefore vanishes, proving that
$q^{(p)}_{M,n}\to0$.

Finally, suppose that $s_{p,n}\to+\infty$.
Then $x/\log n\geq1/(2p)$ for all sufficiently large $n$,
so \eqref{eq:cumulants-spn} with $q=p$ gives
\[
    \frac{K_p}{K_1^p}\longrightarrow\infty.
\]
By Theorem~\ref{thm:radial-dominance},
$K_q\sim A_{qM}(n)>0$ for every $2\leq q\leq p$.
Thus all terms in the moment-cumulant formula are
nonnegative for sufficiently large $n$.
Subtracting the term $m_1=p$ yields
\[
    q^{(p)}_{M,n}
    \geq \frac{K_p}{K_1^p}
    \longrightarrow\infty.
\]
This completes the proof of
\eqref{eq:q-intermediate-powers}.
\end{proof}

\subsection{\texorpdfstring{Small Powers: Proof of \eqref{eq:q-small-powers}}{Small Powers}}
\label{subsec:q-small-powers}

\begin{proof}[Proof of \eqref{eq:q-small-powers}]
   Since $M=O(\sqrt n)$, every subsequence has a further
subsequence along which either $M=o(\sqrt n)$ or
$M\asymp\sqrt n$. It therefore suffices to prove
\eqref{eq:q-small-powers} in these two cases.

  \noindent  \textbf{Case 1:} $M=o(\sqrt{n})$. Here we can use a crude upper bound on the $p$-moment of $Y_{M,n}$. Recall the  expression for $\E [Y^p_{M,n}]$ established in~\eqref{eq:y-perm}.
    Since all the summands are non-negative, we can upper-bound the term in \eqref{eq:y-perm} by the sum over all permutations:
    \begin{align*}
     \E[Y_{M,n}^p]
     &\leq\frac{p!}{n^{pM}}\sum_{\pi\in S_{pM}}n^{\cyc(\pi)}
     =\frac{p!}{n^{pM}}\prod_{j=0}^{pM-1}(n+j)
     \\ &=p!\prod_{j=0}^{pM-1}\left(1+\frac jn\right)
     \leq p!\exp\left(\frac{pM(pM-1)}{2n}\right) \\ & \to p!.
    \end{align*}
    Here, we used the cycle-counting identity (see, e.g., \citep[pp.~206--213]{Comtet:1974})
    \begin{equation}
    \sum_{\pi\in S_{pM}}n^{\cyc(\pi)}=\prod_{j=0}^{pM-1}(n+j)
    \end{equation}
     and the inequality $1+u\leq e^u$, $u\geq 0$. %
    By Proposition~\ref{prop:bridge},
$\mathbb E[Y_{M,n}]\to1$.
Moreover, Theorem~\ref{thm:bridge-variance} and
$R^{(p)}_{k,n}\geq0$ imply $q^{(p)}_{M,n}\geq0$.
Consequently,
\[
    0\leq q^{(p)}_{M,n}
    \leq
    \frac{p!\exp\!\left(pM(pM-1)/(2n)\right)}
         {(\mathbb E[Y_{M,n}])^p}
    -p!
    \longrightarrow0.
\]

  \noindent  \textbf{Case 2:} $M \asymp \sqrt{n}$. 
By Theorem~\ref{thm:radial-dominance}, part~1,
for every fixed $2\leq q\leq p$,
\[
    A_{qM}(n)
    =\frac{n}{qM+1}
      \prod_{j=1}^{qM}\left(1+\frac jn\right)
    \leq
    \frac{n}{qM+1}
    \exp\left(\frac{qM(qM+1)}{2n}\right)
    =O(\sqrt n).
\]
Here we used $n/M=O(\sqrt n)$ and $M^2/n=O(1)$.
Theorem~\ref{thm:radial-dominance}, part~2,
therefore yields
\begin{equation}\label{n1}
    |K_q|\leq(1+C_q)A_{qM}(n)=O(\sqrt n).
\end{equation}

    On the other hand, %
    Proposition~\ref{prop:bridge}
and $M\asymp\sqrt n$ imply
\begin{equation} \label{n2}
    K_1=M\mathbb E[Y_{M,n}]=\Theta(\sqrt n).
\end{equation}
Therefore, by combining \eqref{n1} and \eqref{n2}, we have
    \begin{equation}
        \dfrac{K_q}{K_1^q} = O(n^{(1-q)/2}) \to 0 \quad \text{for all} \quad q\geq 2.
    \end{equation}
In the moment-cumulant expansion
\eqref{eq:moment-cumulant}, the term corresponding to
$m_1=p$ equals $p!$.
Every other term contains at least one factor
$K_q/K_1^q$ with $2\leq q\leq p$ and therefore tends
to zero. Since there are only finitely many terms,
we conclude that $q^{(p)}_{M,n}\to0$.
\end{proof}

\subsection{\texorpdfstring{Large Powers: Proof of \eqref{eq:q-large-powers}}{Large Powers}}
\label{subsec:q-large-powers}

\begin{proof}[Proof of \eqref{eq:q-large-powers}]
Since $M\geq c n^{2/3}$, we have $M/\sqrt n\to\infty$.
By Theorem~\ref{thm:radial-dominance}, part~3,
$K_q\sim A_{qM}(n)>0$ for every $2\leq q\leq p$.
Thus all terms in \eqref{eq:moment-cumulant} are
nonnegative for sufficiently large $n$.
Subtracting the term corresponding to $m_1=p$, which
equals $p!$, gives
\begin{align} \label{eq:Kp_lower_bound}
    q^{(p)}_{M,n}\geq\frac{K_p}{K_1^p},
    \qquad
    K_p\geq\frac12 A_{pM}(n).
\end{align}
For the denominator, we use $K_1=M\E Y_{M,n}$ and the exact mean
formula \eqref{eq:mean_expression}:
\begin{align} \label{eq:K1_upper_bound}
 K_1=\frac n{M+1}\left[
     \prod_{j=0}^{M}\left(1+\frac jn\right)
    -\prod_{j=0}^{M}\left(1-\frac jn\right)\right]
\leq\frac n{M+1}\prod_{j=0}^{M}\left(1+\frac jn\right)
 =\frac{(n+M)!}{(M+1)(n-1)!\,n^M}.
\end{align}
Indeed, the subtracted product is nonnegative when $M<n$ and is zero when $M\geq n$, because it then contains the factor $1-n/n$. Combining \eqref{eq:Kp_lower_bound}, \eqref{eq:K1_upper_bound}, and Theorem~\ref{thm:radial-dominance}, part~1, yields:
\begin{align}
 q^{(p)}_{M,n}
 \geq\frac12\frac{(M+1)^p}{pM+1}
       \frac{(n+pM)!\,((n-1)!)^{p-1}}{((n+M)!)^p}\notag
 =\frac{(M+1)^p}{2n^{p-1}(pM+1)}
 \frac{\displaystyle\prod_{j=1}^{pM}(1+j/n)}
      {\displaystyle\left[\prod_{j=1}^{M}(1+j/n)\right]^p}.
 \label{eq:large-power-mean-ratio}
\end{align}
The product ratio has a simple interpretation: split its numerator into $p$ consecutive blocks of length $M$. The first block cancels, and each later block is larger than the denominator block:
\[
 \frac{\prod_{j=1}^{pM}(1+j/n)}
      {[\prod_{j=1}^{M}(1+j/n)]^p}
 =\prod_{a=1}^{p-1}\prod_{j=1}^{M}
       \left(1+\frac{aM}{n+j}\right).
\]
Using $\log(1+u)\geq u/(1+u)$ and $n+j+aM\leq n+pM$,
\[
 \sum_{a=1}^{p-1}\sum_{j=1}^{M}
 \log\left(1+\frac{aM}{n+j}\right)
 \geq\frac{p(p-1)M^2}{2(n+pM)}.
\]
Also $(M+1)^p/(pM+1)\geq1/p$. Hence, for all sufficiently large $n$,
\begin{equation}\label{eq:large-power-lower-bound}
\log q^{(p)}_{M,n}
\geq\frac{p(p-1)M^2}{2(n+pM)}
                  -(p-1)\log n-\log(2p).
\end{equation}
For $M\geq c n^{2/3}$, the first term is bounded below by a positive constant times $n^{1/3}$, since $M^2/(n+pM)$ increases with $M$. %
Since $n^{1/3}/\log n\to\infty$, we conclude that
$q^{(p)}_{M,n}\to\infty$.
\end{proof}

\subsection{Completion of Theorem~\ref{thm:q-asymptotics}}
Sections~\ref{subsec:q-intermediate-powers}--\ref{subsec:q-large-powers}
establish the asymptotics of $q^{(p)}_{M,n}$ in three power ranges.
We now combine these results, transfer them to $R^{(p)}_{k,n}$
using the bridge, and establish the $L^p$-convergence asserted in
Theorem~\ref{thm:q-asymptotics}.

\noindent\textbf{Combining the power regimes.}
For $M=O(\sqrt n)$, we necessarily have $s_{p,n}\to-\infty$, and~\eqref{eq:q-small-powers} gives $q^{(p)}_{M,n}\to0$.
For $M\geq c n^{2/3}$ with fixed $c>0$, we necessarily have $s_{p,n}\to+\infty$, and~\eqref{eq:q-large-powers} gives
$q^{(p)}_{M,n}\to\infty$.
In the intermediate range $\sqrt n\ll M\ll n^{2/3}$, all three asymptotic behaviors for $s_{p,n}$ are covered by~\eqref{eq:q-intermediate-powers}. In particular, the critical window $s_{p,n}\to s\in\mathbb R$ lies in this range. By a standard subsequence argument, these conclusions also cover sequences that do not stay in a single
power range. 

\noindent\textbf{Transferring the limits via the singular-eigenvalue bridge.}
Theorem~\ref{thm:bridge-variance} gives
\[
    \frac{1}{p^p}q^{(p)}_{M,n}
    \leq R^{(p)}_{k,n}
    \leq \frac{1}{p!}q^{(p)}_{M,n}.
\]
Thus $q^{(p)}_{M,n}\to0$ implies $R^{(p)}_{k,n}\to0$,
and $q^{(p)}_{M,n}\to\infty$ implies $R^{(p)}_{k,n}\to\infty$.
In the critical window, $q^{(p)}_{M,n}$ converges to a strictly
positive finite limit, so the same bounds give $R^{(p)}_{k,n}=\Theta(1)$.
This establishes the excess moment bounds in the theorem.

\noindent\textbf{Convergence in $L^p$.}
It remains to establish concentration when $s_{p,n}\to-\infty$.
Set $X_{k,n}:=Z_{k,n}/\mathbb E Z_{k,n}$, so that
$X_{k,n}\geq0$ and $\mathbb E X_{k,n}=1$.
For every integer $p\geq2$,
\[
    |x-1|^p\leq x^p-1-p(x-1),\qquad x\geq0.
\]
Taking expectations, the linear term vanishes and we obtain
\[
    \mathbb E|X_{k,n}-1|^p
    \leq \mathbb E[X_{k,n}^p]-1
    =R^{(p)}_{k,n}\longrightarrow0.
\]
Hence $Z_{k,n}/\mathbb E Z_{k,n}\to1$ in $L^p$, completing
the proof of Theorem~\ref{thm:q-asymptotics}.

\subsection{Proof of Corollary~\ref{corr:critical_convergence_Zp}}\label{subsec:proof-raw}
\begin{proof}
For $p=1$, the claims follow from Proposition~\ref{prop:bridge}. For $p\geq2$, use
\[
 \mathbb E[Z_{k,n}^p]=(\mathbb E Z_{k,n})^p(1+R^{(p)}_{k,n}).
\]
When $k=O(\sqrt n)$, Theorem~\ref{thm:q-asymptotics} gives $R^{(p)}_{k,n}\to0$, and Proposition~\ref{prop:bridge} gives the stated limits. When $k/\sqrt n\to\infty$, Jensen's inequality gives $\mathbb E[Z_{k,n}^p]\geq(\mathbb E Z_{k,n})^p\to\infty$. Finally, if $s_{p,n}\to-\infty$ and $k\gg\sqrt n$, the same identity and the intermediate-power mean asymptotic give the claimed equivalence.
\end{proof}

\section{Proof of Radial Dominance} \label{sec:proof_thm_radial}
Recall the definition of $A_{pM}(n)$ in \eqref{eq:def_A}.
The proof of Theorem \ref{thm:radial-dominance} consists of the following three lemmas, which show that the mixed cumulant $K_p$ defined in \eqref{eq:def_K} contains the radial term $A_{pM}(n)$ as one of its summands (Lemma~\ref{lem:cumulant-trace-expansion}), derive the scaling of this term (Lemma~\ref{lem:radial-term}), and show that all the remaining terms scale not faster than $A_{pM}(n)$ (Lemma~\ref{lem:nonradial-trace-bound}). Each lemma is proved immediately after its statement. We then assemble the theorem using these lemmas in Section~\ref{sec:proof-radial-dominance-completion}.

\subsection{Radial Term Expression and Scaling}
\label{subsec:proof-radial-term}

\begin{lemma}[Radial term]\label{lem:radial-term}
For every integer $L\geq0$,
\begin{equation}
\label{eq:radial-term-exact-detail}
    A_L(n)
=
    \frac{1}{n^L}
    \sum_{j=0}^{n-1}\frac{(j+L)!}{j!}
    =
    \frac{(n+L)!}{(L+1)(n-1)!\,n^L}.
\end{equation}
In particular, for every fixed $p\geq1$, if $M=o(n^{2/3})$ and $M\to\infty$, then
\begin{equation}
\label{eq:radial-term-asymptotic}
   A_{pM}(n) \sim
    \frac{n}{pM}
    \exp\left\{
        \frac{p^2M^2}{2n}
    \right\}.
\end{equation}
\end{lemma}

\begin{proof}
The proof amounts to integrating the well-known Ginibre one-point eigenvalue density~\citep{ginibre1965statistical}, given by
\begin{equation}\label{eq:ginibre_one_point_density}
\rho_n^{(1)}(z)=\frac n\pi e^{-n|z|^2}
                   \sum_{j=0}^{n-1}\frac{(n|z|^2)^j}{j!}.
\end{equation}
Here the normalization is such that the density integrates to $n$, so the defining identity is
$\E\sum_i f(\lambda_i)=\int_{\mathbb C}f(z)\rho_n^{(1)}(z)\,d^2z$. Applying this identity with $f(z)=|z|^{2L}$, polar coordinates and
the substitution $t=nr^2$ give
\begin{align*}
 A_L(n)
 &=2n\sum_{j=0}^{n-1}\frac{n^j}{j!}
       \int_0^\infty r^{2(L+j)+1}e^{-nr^2}\,dr
       \\
 &=\frac1{n^L}\sum_{j=0}^{n-1}\frac1{j!}
       \int_0^\infty t^{L+j}e^{-t}\,dt
 =\frac1{n^L}\sum_{j=0}^{n-1}\frac{(L+j)!}{j!}.
\end{align*}
Then, we can get the expression for $A_L(n)$ in \eqref{eq:radial-term-exact} using the hockey-stick identity as follows:
\begin{equation}
     \sum_{j=0}^{n-1}\frac{(j+L)!}{j!}
 =L!\sum_{j=0}^{n-1}\binom{j+L}{j}
 =L!\binom{n+L}{n-1}
 =\frac{(n+L)!}{(L+1)(n-1)!} = \dfrac{n^{L+1}}{(L+1)}\prod_{j=1}^L\left(1 + \frac{j}{n}\right).
\end{equation}
It remains to prove \eqref{eq:radial-term-asymptotic}.
For $L=pM=o(n^{2/3})$, the expansion
$\log(1+t)=t+O(t^2)$ yields
\begin{equation}
 \log\prod_{j=1}^{L}\left(1+\frac jn\right)
 =\frac{L(L+1)}{2n}+O\left(\frac{L^3}{n^2}\right)
 =\frac{L^2}{2n}+o(1).
\end{equation}
This gives the desired asymptotic in \eqref{eq:radial-term-asymptotic}:
\begin{equation}
    A_{pM}(n) = \dfrac{n}{pM+1}\prod_{j=1}^{pM}\left(1 + \frac{j}{n}\right) \sim \dfrac{n}{pM} \exp\left(\dfrac{p^2M^2}{2n}\right).
\end{equation}
\end{proof}

\subsection{Mixed Cumulant Expansion}
\begin{lemma}[Trace expansion of the mixed cumulant]
\label{lem:cumulant-trace-expansion}
For nonnegative integers $a,b$, define the $n\times n$ matrix $B_{a,b}= B_{a,b}(n,M)$ by
\begin{equation}
\label{eq:B-matrices}
    (B_{a,b})_{ij}
    :=
    \mathbf{1}_{\{j + bM=i+aM\}}
    \frac{(i+aM)!}
    {\sqrt{i!j!}\,n^{(a+b)M/2}},
    \qquad 0\leq i,j\leq n-1.
\end{equation}
In particular, $B_{0,0}=I$. For every fixed $p\geq1$,
\begin{equation}\label{eq:cumulant_expansion}
K_p
=
(p!)^2
\sum_{\ell=1}^{2p}
\frac{(-1)^{\ell-1}}{\ell}
\sum_{\substack{
a_1+\cdots+a_\ell=p\\
b_1+\cdots+b_\ell=p\\
a_r+b_r>0}}
\frac{
\tr(B_{a_1,b_1}\cdots B_{a_\ell,b_\ell})
}{
\prod_{r=1}^\ell a_r!b_r!
}.
\end{equation}
Moreover, the term corresponding to $\ell=1$ is exactly $A_{pM}(n)$.
\end{lemma}

\begin{proof}
All generating functions in this proof are understood as formal power series, whose coefficients are finite moments. Recall that the joint cumulants of $T=\tr(\W^M)$ and $\overline T$ are defined by the coefficients of the cumulant generating function:
\begin{equation} \label{eq:cumulant_expansion-detail}
     \log\E e^{sT+t\overline T}
 =\sum_{\substack{a,b\geq0\\a+b>0}}
 \frac{s^at^b}{a!b!}
 \kappa(\underbrace{T,\ldots,T}_{a},
        \underbrace{\overline T,\ldots,\overline T}_{b}) = \sum_{q\geq1}\frac{K_q}{(q!)^2}s^qt^q.
\end{equation}
Here all the terms with $a\neq b$ are zero by rotational invariance. Thus, for the proof, it suffices to compute $\E e^{sT+t\overline T}$, take its logarithm, and extract the coefficient of $s^pt^p$. 

Since $T = \sum_{i=1}^n \lambda_i^M$, we can compute the expectation $\E e^{sT+t\overline T}$ using the joint eigenvalue density of the complex Ginibre ensemble~\citep{ginibre1965statistical}. We write the density using the Vandermonde determinant
and normalized monomials as follows:
\begin{equation*}
    \rho_n(z_1,\ldots,z_n)
 =\frac1{n!}
   \left|\det[\phi_{i-1}(z_j)]_{i,j=1}^n\right|^2
   \prod_{j=1}^n w_n(z_j), \quad 
 w_n(z)=\frac n\pi e^{-n|z|^2},
 \quad
 \phi_i(z)=\frac{n^{i/2}z^i}{\sqrt{i!}}.
\end{equation*}
Since $\rho_n$ is the density, it holds for every symmetric integrable function $F:\mathbb{C}^n \to \mathbb{C}$ that:
\begin{equation}\label{eq:joint-density-expectation}
 \E F(\lambda_1,\ldots,\lambda_n)
 =\int_{\mathbb C^n}F(z_1,\ldots,z_n)
       \rho_n(z_1,\ldots,z_n)\prod_{j=1}^n d^2z_j.
\end{equation}
Let $f_i,h_i,g:\mathbb C\to\mathbb C$ ($1\leq i\leq n$) be measurable functions satisfying
\begin{align*}
    \int_{\mathbb C}
    |f_i(z)h_j(z)g(z)|\,w_n(z)\,d^2z<\infty,
    \qquad 1\leq i,j\leq n.
\end{align*} 
Then, Andr\'{e}ief's identity~\citep{andreief:1886,forrester2019meet} says
\begin{equation} \label{eq:left_hand}
 \frac1{n!}\int_{\mathbb C^n}
   \det[f_i(z_j)]\det[h_i(z_j)]
   \prod_{j=1}^n g(z_j)\,w_n(z_j)\,d^2z_j
 =\det[C_{ij}]_{i,j=1}^n,
\end{equation}
where
\begin{equation}
  C_{ij}=\int_{\mathbb C}f_i(z)h_j(z)g(z)\,w_n(z)\,d^2z.
\end{equation}
Now take $f_i=\phi_{i-1}$, $h_i=\overline{\phi_{i-1}}$, and
$g(z)=e^{sz^M+t\bar z^M}$. We apply \eqref{eq:left_hand} to
$g$ coefficientwise in $s$ and $t$.
Each coefficient of $f_i(z)h_j(z)g(z)$ is a polynomial
in $z$ and $\bar z$ and is therefore absolutely integrable
against $w_n(z)\,d^2z$.
Accordingly, $\mathbb E e^{sT+t\overline T}$ and its
logarithm are understood as formal power series.

Since
$e^{sT+t\overline T}=\prod_{j=1}^n g(z_j)$, the left-hand side of \eqref{eq:left_hand} is
exactly the expectation under the Ginibre eigenvalue density. We obtain
\begin{equation}\label{eq:formal-det}
 \E e^{sT+t\overline T}
 =\det\left[\int \phi_i(z)\overline{\phi_j(z)}
                         e^{sz^M+t\bar z^M}\,w_n(z)\,d^2z\right]_{i,j=0}^{n-1}.
\end{equation}
Now we want to write the matrix in \eqref{eq:formal-det} as $I+X$, where
\[
 X=\sum_{a+b>0}\frac{s^at^b}{a!b!}B_{a,b}.
\]
Here the matrices $B_{a,b}$ are the ones given in the lemma statement. They arise from~\eqref{eq:formal-det} as follows:
\begin{align*}
(B_{a,b})_{ij} = \int_{\mathbb C}\phi_i(z)\overline{\phi_j(z)}
                   z^{aM}\bar z^{bM}\,w_n(z)\,d^2z
 &=\frac{n^{(i+j)/2}}{\sqrt{i!j!}}
   \int_{\mathbb C}z^{i+aM}\bar z^{j+bM}\,w_n(z)\,d^2z\\
 &=\mathbf1_{\{i+aM=j+bM\}}
   \frac{(i+aM)!}{\sqrt{i!j!}\,n^{(a+b)M/2}},
\end{align*}
where we used the identity
\[
 \int_{\mathbb C}z^u\bar z^v\,w_n(z)\,d^2z
 =\mathbf1_{\{u=v\}}\frac{u!}{n^u},
 \qquad u,v\geq0,
\]
and that $(i+j)/2-(i+aM)=-(a+b)M/2$ whenever the indicator
is nonzero. In particular, $B_{0,0}=I$.

Taking the formal logarithm of the expression in~\eqref{eq:formal-det} gives
\begin{align} \label{eq:trace_expansion}
 \log\E e^{sT+t\overline T}
 =\log\det(I+X)
 =\sum_{\ell\geq1}\frac{(-1)^{\ell-1}}{\ell}\tr(X^\ell).
\end{align}
We now compare the coefficients of $s^pt^p$ in the cumulant
generating function \eqref{eq:cumulant_expansion-detail} and its trace expansion \eqref{eq:trace_expansion}. The former coefficient
is $K_p/(p!)^2$. On the other hand, every monomial in $s$ and $t$
in $X^\ell$ has total degree at least $\ell$, so only $\ell\leq2p$
can contribute. Expanding the product $X^\ell$ and multiplying the
resulting coefficient by $(p!)^2$ gives
\[
 K_p=(p!)^2
 \sum_{\ell=1}^{2p}\frac{(-1)^{\ell-1}}{\ell}
 \sum_{\substack{a_1+\cdots+a_\ell=p\\
                  b_1+\cdots+b_\ell=p\\
                  a_r,b_r\geq0,\ a_r+b_r>0}}
 \frac{\tr(B_{a_1,b_1}\cdots B_{a_\ell,b_\ell})}
      {\prod_{r=1}^{\ell}a_r!b_r!},
\]
which is exactly the desired expression~\eqref{eq:cumulant_expansion}.
For $\ell=1$, the only pair is $(a_1,b_1)=(p,p)$, and its coefficient
is $(p!)^2/(p!p!)=1$. Its trace is
\[
 \tr B_{p,p}=\frac1{n^{pM}}\sum_{j=0}^{n-1}\frac{(j+pM)!}{j!}
 =A_{pM}(n),
\]
by Lemma~\ref{lem:radial-term}, which completes the proof.
\end{proof}

\subsection{Upper Bound on Non-Radial Terms} \label{subsec:proof-nonradial-trace-bound}
\begin{lemma}[Non-radial trace bound]
\label{lem:nonradial-trace-bound}
Let $p\geq 2$ and $\ell\geq2$, and let
$a_1,\ldots,a_\ell,b_1,\ldots,b_\ell$ be nonnegative integers such that $\sum_{r=1}^\ell a_r = \sum_{r=1}^\ell b_r = p$ and $a_r+b_r>0$ for all $r$.
Then for the matrices $B$ defined in~\eqref{eq:B-matrices} and any finite $n, M$:
\label{eq:nonradial-trace-bound}
\begin{equation}\label{eq:trace-bound}
 0\leq\tr(B_{a_1,b_1}\cdots B_{a_\ell,b_\ell})
 \leq A_{pM}(n)\exp\left(-\frac{(p-1)M^2}{n+pM}\right).
\end{equation}
\end{lemma}

\begin{proof}
We begin by expanding the trace as follows:
\begin{equation}
\tr(B_{a_1,b_1}\cdots B_{a_\ell,b_\ell})
 \quad=\sum_{i_1,\ldots,i_\ell=0}^{n-1}
 (B_{a_1,b_1})_{i_1i_2}
 (B_{a_2,b_2})_{i_2i_3}\cdots
 (B_{a_\ell,b_\ell})_{i_\ell i_1}.
\end{equation}
By definition of the matrices $B_{a,b}$ in ~\eqref{eq:B-matrices}, the lower bound follows immediately. Moreover, by~\eqref{eq:B-matrices}, the term with indices $(i_1,\dots, i_\ell)$ is non-zero only when the indices satisfy 
\begin{equation} \label{eq:def_i_r}
     i_{r+1}+b_rM=i_r+a_rM,
 \qquad\text{equivalently}\qquad
 i_{r+1}=i_r+(a_r-b_r)M,
 \quad 1\leq r\leq\ell.
\end{equation}
Here we use the convention $i_{\ell+1}=i_1$, which is consistent with the above equations because $\sum a_r = \sum b_r = p$. For each admissible index tuple $(i_1,\ldots,i_\ell)$,
the corresponding summand in the trace expansion equals
\begin{equation}\label{eq:equality_trace_summand}
\begin{aligned}
 &(B_{a_1,b_1})_{i_1i_2}\cdots(B_{a_\ell,b_\ell})_{i_\ell i_1}\\
 &\qquad=\prod_{r=1}^{\ell}
 \frac{(i_r+a_rM)!}{\sqrt{i_r!i_{r+1}!}\,n^{(a_r+b_r)M/2}}
 =\frac1{n^{pM}}\prod_{r=1}^{\ell}\frac{(i_r+a_rM)!}{i_r!}.
\end{aligned}
\end{equation}
where we used $\prod_r\sqrt{i_r!i_{r+1}!}=\prod_r i_r!$, and $\sum_r(a_r+b_r)M/2=pM$. 
We can now associate every non-zero term to its maximal index, defined as follows:
\begin{equation}
    j := \max_{1\leq r\leq\ell}i_r 
    = i_1 +M \max_{1\leq r\leq\ell} \sum_{s=1}^{r-1}(a_s - b_s),
\end{equation}
where the second equality follows by iteratively using \eqref{eq:def_i_r}. 
Clearly, $j$ uniquely determines $i_1$ since $a_s,b_s$ are fixed. Moreover, notice that $j$ can be analogously rewritten via any index $i_k$ as $$j= i_k +  M \max_{1 \le r \le \ell}\sum_{s=k}^{k+r-2}\left(a_{((s-1)\bmod \ell)+1}- b_{((s-1)\bmod \ell)+1} \right).$$ Therefore, $j$ uniquely determines the whole tuple of indices $(i_1, \dots, i_\ell)$. Let $\mathcal B_j$ denote the contribution of this tuple, and set $\mathcal B_j=0$ if no admissible tuple exists for a given $j$. Then, using \eqref{eq:equality_trace_summand}, we may write
\begin{equation}
    \tr(B_{a_1,b_1}\cdots B_{a_\ell,b_\ell}) = \sum_{j=0}^{n-1} \mathcal{B}_j, 
    \qquad \mathcal B_j
=\frac1{n^{pM}}\prod_{r=1}^{\ell}\frac{(i_r+a_rM)!}{i_r!}
\leq\frac1{n^{pM}}\prod_{r=1}^{\ell}\frac{(j+a_rM)!}{j!}.
\end{equation}
Recall that the radial term $A_{pM}(n)$ consists of similar terms by Lemma~\ref{lem:radial-term}. Namely, we can write:
\begin{equation}
    A_{pM}(n) = \sum_{j=0}^{n-1}\mathcal{A}_j, \qquad \mathcal{A}_j := \dfrac{1}{n^{pM}}\dfrac{(j+pM)!}{j!}.
\end{equation}
Thus, in what follows, we will bound the product trace terms via the radial terms as follows:
\begin{equation} \label{eq_claim}
    \mathcal{B}_j \leq \mathcal{A}_j \exp\left(-\frac{(p-1)M^2}{n+pM}\right), \qquad j=0,\dots, n-1.
\end{equation}
This is immediate for $j$ such that $\mathcal{B}_j=0$. Otherwise, there are two possible cases. 

\emph{Case 1: at least two of the $a_r$ are positive.}
Among such configurations with $\sum_ra_r=p$, the product
$\prod_r(j+a_rM)!/j!$ is largest at $(p-1,1,0,\ldots,0)$,
up to permutation. Therefore,
\begin{align*}
 \frac{\mathcal B_j}{\mathcal A_j}
 \leq\frac{(j+(p-1)M)!(j+M)!}{j!(j+pM)!}
 =\prod_{v=1}^{M}
 \left(1-\frac{(p-1)M}{j+(p-1)M+v}\right)
 \leq\exp\left(-\frac{(p-1)M^2}{n+pM}\right).
\end{align*}
The last step uses $1-u\leq e^{-u}$ and
$j+(p-1)M+v\leq n+pM$. This proves
the claim \eqref{eq_claim} in the first case.

\emph{Case 2: exactly one $a_r$ is positive.} In this case, $a_r=p$ for some $r$ and all the other factorial ratios in $\mathcal B_j$ equal one. Since $\ell\geq 2$, we must have $b_r < p$. Otherwise, all other pairs would be $(0,0)$, which corresponds to the $\ell=1$ term. Then we have the following for the corresponding tuple index $i_r$:
\begin{equation}
i_{r+1}=i_r+(p-b_r)M\geq i_r+M,
 \qquad\text{so}\qquad i_r\leq j-M.
\end{equation}
This allows to derive an upper bound on $\mathcal{B}_j$:
\begin{equation}
    \mathcal{B}_j = \dfrac{1}{n^{pM}}\dfrac{(i_r+pM)!}{i_r!} \leq \dfrac{1}{n^{pM}}\dfrac{(j - M +pM)!}{(j-M)!}.
\end{equation}
Finally, analogous to the previous case, this gives the desired upper bound on the summands ratio:
\begin{equation}
    \frac{\mathcal B_j}{\mathcal A_j}
 \leq \prod_{v=1}^{M}\left(1-\frac{pM}{j+(p-1)M+v}\right)
 \leq\exp\left(-\frac{pM^2}{n+pM}\right)
 \leq\exp\left(-\frac{(p-1)M^2}{n+pM}\right).
\end{equation}

Thus, the upper bound \eqref{eq_claim} holds in both cases, and summing over $j$ completes the proof:
\begin{equation}
 \tr(B_{a_1,b_1}\cdots B_{a_\ell,b_\ell})
 =\sum_{j=0}^{n-1}\mathcal B_j
 \leq e^{-\frac{(p-1)M^2}{(n+pM)}}\sum_{j=0}^{n-1}\mathcal A_j
 =A_{pM}(n)e^{-\frac{(p-1)M^2}{(n+pM)}}.
\end{equation}

\end{proof}

\subsection{Completion of the Radial-Dominance Theorem}\label{sec:proof-radial-dominance-completion}
\begin{proof}[Proof of Theorem \ref{thm:radial-dominance}]
Part~1 follows directly from Lemma~\ref{lem:radial-term}.

By Lemma~\ref{lem:cumulant-trace-expansion}, the radial
term in \eqref{eq:cumulant_expansion} (for $\ell=1$) is exactly
$A_{qM}(n)$. By Lemma~\ref{lem:nonradial-trace-bound},
every non-radial trace in this expansion is bounded by
\[
    A_{qM}(n)
    \exp\left(-\frac{(q-1)M^2}{n+qM}\right).
\]
For fixed $q$, the expansion contains only finitely many
terms, whose coefficients depend only on $q$.
The triangle inequality therefore gives
\[
    |K_q-A_{qM}(n)|
    \leq C_q A_{qM}(n)
    \exp\left(-\frac{(q-1)M^2}{n+qM}\right)
\]
for a constant $C_q<\infty$, proving part~2.

Finally, we turn to the proof of part 3. If $M/\sqrt n\to\infty$, then
\[
    \frac{n+qM}{M^2}
    =\frac{n}{M^2}+\frac{q}{M}
    \longrightarrow0.
\]
Thus the exponential factor in \eqref{eq:radial-error-bound} tends to zero,
and hence $K_q\sim A_{qM}(n)$.
If, in addition, $M=o(n^{2/3})$, then
Lemma~\ref{lem:radial-term} yields
\[
    K_q\sim A_{qM}(n)
    \sim \frac{n}{qM}
    \exp\left(\frac{q^2M^2}{2n}\right).
\]
This completes the proof of part 3.
\end{proof}

\backmatter

\section*{Acknowledgments}
MS is funded by the Deutsche Forschungsgemeinschaft (DFG, German Research Foundation) --- Project number 572957656 (“RMT4DL: Random Matrix Theory for Deep Learning.”). Numerical experiments for this research were conducted on the Lichtenberg high-performance computing cluster of TU Darmstadt. The work of ND is partially supported by  the Aarhus University Research
Foundation (AUFF),  project numbers 47221 and 47388, and by a RECRUIT Sabbatical
Grant from the Novo Nordisk Foundation, project number 52629. ND gratefully acknowledges the hospitality of the Department of
Statistics at Harvard University during the preparation of this
manuscript. 

\section*{Declaration of AI Use}
The authors made use of ChatGPT (GPT-5, GPT-6) in the course of this project. In particular, AI tools were used to edit and proofread the paper text, edit its \LaTeX{} source, assist in writing code for numerical experiments, and help verify, format, and develop proofs. We did \emph{not} use generative AI to develop the original research idea or the main proof strategy, nor did we rely on it for the main literature review beyond assisting with novelty verification. We reviewed all AI-assisted work and take responsibility for the final content of this paper.

\bibliography{references}

@article{seleznova:2026,
  author  = {Seleznova, Mariia},
  title   = {How Long Does Infinite Width Last? Signal Propagation in Long-Range Linear Recurrences},
  journal = {arXiv preprint arXiv:2605.05113},
  year    = {2026}
}

@book{Comtet:1974,
  author    = {Comtet, Louis},
  title     = {Advanced Combinatorics},
  publisher = {D. Reidel Publishing Company},
  address   = {Dordrecht},
  year      = {1974}
}

@article{dubach:2024,
  author  = {Dubach, Guillaume},
  title   = {On the Number of Cycles in Commutators of Random Permutations},
  journal = {The Annals of Applied Probability},
  volume  = {34},
  number  = {4},
  pages   = {4072--4084},
  year    = {2024},
  doi     = {10.1214/24-AAP2059}
}

@book{bai:silverstein:2010,
  author    = {Bai, Zhidong and Silverstein, Jack W.},
  title     = {Spectral Analysis of Large Dimensional Random Matrices},
  series    = {Springer Series in Statistics},
  edition   = {2},
  publisher = {Springer},
  address   = {New York},
  year      = {2010},
  doi       = {10.1007/978-1-4419-0661-8},
  isbn      = {978-1-4419-0660-1}
}

@article{dubach:2021,
  author  = {Dubach, Guillaume and Peled, Yuval},
  title   = {On Words of Non-Hermitian Random Matrices},
  journal = {The Annals of Probability},
  volume  = {49},
  number  = {4},
  pages   = {1886--1916},
  year    = {2021}
}

@article{dubach:2018,
  author  = {Dubach, Guillaume},
  title   = {Powers of Ginibre Eigenvalues},
  journal = {Electronic Journal of Probability},
  volume  = {23},
  pages   = {1--31},
  year    = {2018},
  doi     = {10.1214/18-EJP234}
}

@article{alexeev:2010,
  author  = {Alexeev, Nikita and G{\"o}tze, Friedrich and Tikhomirov, Alexander},
  title   = {Asymptotic Distribution of Singular Values of Powers of Random Matrices},
  journal = {Lithuanian Mathematical Journal},
  volume  = {50},
  number  = {2},
  pages   = {121--132},
  year    = {2010}
}

@article{schoenholz2017deep,
  title={Deep information propagation},
  author={Schoenholz, Samuel S and Gilmer, Justin and Ganguli, Surya and Sohl-Dickstein, Jascha},
  journal={International Conference on Learning Representations},
  year={2017}
}

@inproceedings{chen2018dynamical,
  title={Dynamical isometry and a mean field theory of RNNs: Gating enables signal propagation in recurrent neural networks},
  author={Chen, Minmin and Pennington, Jeffrey and Schoenholz, Samuel},
  booktitle={International Conference on Machine Learning},
  pages={873--882},
  year={2018},
  organization={PMLR}
}

@article{yang2019wide,
	title={Wide feedforward or recurrent neural networks of any architecture are gaussian processes},
	author={Yang, Greg},
	journal={Advances in neural information processing systems},
	volume={32},
	year={2019}
}

@inproceedings{orvieto2023resurrecting,
  title={Resurrecting recurrent neural networks for long sequences},
  author={Orvieto, Antonio and Smith, Samuel L and Gu, Albert and Fernando, Anushan and Gulcehre, Caglar and Pascanu, Razvan and De, Soham},
  booktitle={International Conference on Machine Learning},
  pages={26670--26698},
  year={2023},
  organization={PMLR}
}

@article{gu2021efficiently,
  title={Efficiently modeling long sequences with structured state spaces},
  author={Gu, Albert and Goel, Karan and R{\'e}, Christopher},
  journal={arXiv preprint arXiv:2111.00396},
  year={2021}
}

@article{gu2023mamba,
  title={Mamba: Linear-time sequence modeling with selective state spaces},
  author={Gu, Albert and Dao, Tri},
  journal={arXiv preprint arXiv:2312.00752},
  year={2023}
}

@book{trefethen2005spectra,
  author    = {Lloyd N. Trefethen and Mark Embree},
  title     = {Spectra and Pseudospectra: The Behavior of Nonnormal Matrices and Operators},
  publisher = {Princeton University Press},
  address   = {Princeton, NJ},
  year      = {2005}
}

@inproceedings{glorot2010understanding,
  title={Understanding the difficulty of training deep feedforward neural networks},
  author={Glorot, Xavier and Bengio, Yoshua},
  booktitle={Proceedings of the thirteenth international conference on artificial intelligence and statistics},
  pages={249--256},
  year={2010},
  organization={JMLR Workshop and Conference Proceedings}
}

@article{rider2003limit,
  title={A limit theorem at the edge of a non-Hermitian random matrix ensemble},
  author={Rider, Brian},
  journal={Journal of Physics A: Mathematical and General},
  volume={36},
  number={12},
  pages={3401--3409},
  year={2003}
}

@article{rider2014extremal,
  title={Extremal laws for the real Ginibre ensemble},
  author={Rider, Brian and Sinclair, Christopher D},
  year={2014}
}

@article{ginibre1965statistical,
  title={Statistical ensembles of complex, quaternion, and real matrices},
  author={Ginibre, Jean},
  journal={Journal of Mathematical Physics},
  volume={6},
  number={3},
  pages={440--449},
  year={1965},
  publisher={American Institute of Physics}
}

@book{byun2025progress,
  title={Progress on the Study of the Ginibre Ensembles},
  author={Byun, Sung-Soo and Forrester, Peter J},
  year={2025},
  publisher={Springer},
  address={Singapore}
}

@article{trefethen1997pseudospectra,
  title={Pseudospectra of linear operators},
  author={Trefethen, Lloyd N},
  journal={SIAM review},
  volume={39},
  number={3},
  pages={383--406},
  year={1997},
  publisher={SIAM}
}

@article{chalker1998eigenvector,
  title={Eigenvector statistics in non-Hermitian random matrix ensembles},
  author={Chalker, John T and Mehlig, Bernhard},
  journal={Physical review letters},
  volume={81},
  number={16},
  pages={3367},
  year={1998},
  publisher={APS}
}

@article{bar2026revisiting,
  title={Revisiting Glorot Initialization for Long-Range Linear Recurrences},
  author={Bar, Noga and Seleznova, Mariia and Kutyniok, Gitta and Giryes, Raja and others},
  journal={Advances in Neural Information Processing Systems},
  volume={38},
  pages={55397--55422},
  year={2026}
}

@article{liu2023lyapunov,
  title={Lyapunov Exponent, Universality and Phase Transition for Products of Random Matrices},
  author={Liu, Dang-Zheng and Wang, Dong and Wang, Yanhui},
  journal={Communications in Mathematical Physics},
  volume={399},
  number={3},
  pages={1811--1855},
  year={2023},
  publisher={Springer}
}

@inproceedings{liu2024phase,
  title={Phase transitions for infinite products of large non-Hermitian random matrices},
  author={Liu, Dang-Zheng and Wang, Yanhui},
  booktitle={Annales de l'Institut Henri Poincare (B) Probabilites et statistiques},
  volume={60},
  number={4},
  pages={2813--2848},
  year={2024},
  organization={Institut Henri Poincar{\'e}}
}

@article{akemann2020universality,
  title={Universality of local spectral statistics of products of random matrices},
  author={Akemann, Gernot and Burda, Zdzislaw and Kieburg, Mario},
  journal={Physical Review E},
  volume={102},
  number={5},
  pages={052134},
  year={2020},
  publisher={APS}
}

@article{hanin2020products,
  title={Products of Many Large Random Matrices and Gradients in Deep Neural Networks: B. Hanin, M. Nica},
  author={Hanin, Boris and Nica, Mihai},
  journal={Communications in Mathematical Physics},
  volume={376},
  number={1},
  pages={287--322},
  year={2020},
  publisher={Springer}
}

@article{hanin2018neural,
  title={Which neural net architectures give rise to exploding and vanishing gradients?},
  author={Hanin, Boris},
  journal={Advances in neural information processing systems},
  volume={31},
  year={2018}
}

@article{hanin2019finite,
  title={Finite depth and width corrections to the neural tangent kernel},
  author={Hanin, Boris and Nica, Mihai},
  journal={International Conference on Learning Representations},
  year={2020}
}

@inproceedings{seleznova2022neural,
  title={Neural tangent kernel beyond the infinite-width limit: Effects of depth and initialization},
  author={Seleznova, Mariia and Kutyniok, Gitta},
  booktitle={International Conference on Machine Learning},
  pages={19522--19560},
  year={2022},
  organization={PMLR}
}

@article{bourgade2020distribution,
  title={The distribution of overlaps between eigenvectors of Ginibre matrices},
  author={Bourgade, Paul and Dubach, Guillaume},
  journal={Probability Theory and Related Fields},
  volume={177},
  number={1},
  pages={397--464},
  year={2020},
  publisher={Springer}
}

@article{egorychev1981proof,
  title={Proof of the van der Waerden conjecture for permanents},
  author={Egorychev, G. P.},
  journal={Siberian Mathematical Journal},
  volume={22},
  number={6},
  pages={854--859},
  year={1981},
  publisher={Springer}
}

@article{isserlis1918formula,
	title={On a formula for the product-moment coefficient of any order of a normal frequency distribution in any number of variables},
	author={Isserlis, Leon},
	journal={Biometrika},
	volume={12},
	number={1/2},
	pages={134--139},
	year={1918},
	publisher={JSTOR}
}

@article{andreief:1886,
  author  = {Andr{\'e}ief, C.},
  title   = {Note sur une relation entre les int{\'e}grales d{\'e}finies des produits des fonctions},
  journal = {M{\'e}moires de la Soci{\'e}t{\'e} des Sciences Physiques et Naturelles de Bordeaux},
  series  = {3},
  volume  = {2},
  pages   = {1--14},
  year    = {1886}
}

@article{forrester2019meet,
  title={Meet andr{\'e}ief, bordeaux 1886, and andreev, kharkov 1882--1883},
  author={Forrester, Peter J},
  journal={Random Matrices: Theory and Applications},
  volume={8},
  number={02},
  pages={1930001},
  year={2019},
  publisher={World Scientific}
}

\end{document}